\documentclass[twoside,11pt]{article}
\usepackage{jmlr2e_arxiv}

\usepackage[T1]{fontenc}
\usepackage{amsmath}
\usepackage{amssymb}
\usepackage{faktor} 
\usepackage{float}
\usepackage{placeins}
\usepackage[inline]{enumitem}
\usepackage{longtable}
\usepackage{booktabs}
\usepackage{microtype}
\usepackage{tikz-cd}
\usepackage{siunitx}
\usepackage{booktabs}
\usepackage{subcaption}

\usepackage{graphicx}
\usepackage{pgfplots}
\pgfplotsset{compat=1.8}
\usetikzlibrary{math,calc}

\usepackage{breakurl}
\usepackage{mathtools}
\usepackage[utf8]{inputenc}

\usepackage{thmtools}

\newcommand{\pathspVar}[1]{\mathcal{C}^{#1}}

\newcommand{\pathsreduced}[1]{\mathcal{P}^{#1}}

\newcommand{\TSeries}[1]{\prod_{m\geq 0}{#1}^{\otimes m}}
\newcommand{\TSeriesTruncated}[2]{\prod_{m=0}^{#2}{#1}^{\otimes m}}
\newcommand{\FSeries}[1]{\mathbf{T}({#1})}
\newcommand{\FSeriesOne}[1]{\mathbf{T}_{1}({#1})}

\newcommand{\Ker}{\operatorname{Ker}}

\newcommand{\id}{\operatorname{id}}

\newcommand{\linear}{{\operatorname{linear}}}

\newcommand{\calF}{\mathcal{F}}
\newcommand{\calG}{\mathcal{G}}
\newcommand{\calH}{\mathcal{H}}

\newcommand{\calM}{\mathcal{M}}

\newcommand{\calP}{\mathcal{P}}

\newcommand{\calS}{\mathcal{S}}

\newcommand{\calX}{\mathcal{X}}
\newcommand{\calY}{\mathcal{Y}}

\newcommand{\frakT}{\mathfrak{T}}

\newcommand{\RR}{\ensuremath{\mathbb{R}}}
\newcommand{\PP}{\ensuremath{\mathbb{P}}}

\newcommand{\EE}{\ensuremath{\mathbb{E}}}

\newcommand{\bft}{\mathbf{t}}
\newcommand{\bfs}{\mathbf{s}}
\newcommand{\bfx}{\mathbf{x}}
\newcommand{\bfy}{\mathbf{y}}
\newcommand{\bfone}{\mathbf{1}}

\newcommand{\Var}{\operatorname{--var}}

\newcommand{\Sig}{\operatorname{S}}

\newcommand{\kernel}{{\operatorname{k}}} 
\newcommand{\kernelTrunc}[1]{\operatorname{k}_{#1}} 

\newcommand{\floor}[1]{\lfloor #1 \rfloor}

\newcommand{\scal}[1]{\langle #1 \rangle}

\newcommand{\Normal}{\Lambda}

\newcommand{\Norm}[1]{\left\lVert#1\right\rVert}
\newcommand{\norm}[1]{\left\lvert#1\right\rvert}
\newcommand{\bX}{\mathbf{X}}
\newcommand{\bY}{\mathbf{Y}}
\newcommand{\bZ}{\mathbf{Z}}

\newcommand{\IF}{\operatorname{IF}}

\usepackage{lmodern}
\usepackage{todonotes}
\usepackage{booktabs}
\usepackage{float}
\usepackage{bookmark}
\usepackage{appendix}
\newcommand{\ubold}[1]{\fontseries{b}\selectfont#1}
\newcommand\tablescale{0.6}

\usepackage{lastpage}

\firstpageno{1}

\begin{document}




\title{Signature Moments to Characterize Laws of Stochastic Processes}

\author{\name Ilya Chevyrev \email ichevyrev@gmail.com \\
       \addr School of Mathematics \\
       University of Edinburgh\\
       Edinburgh, EH9 3FD, United Kingdom 
       \AND
       \name Harald Oberhauser \email oberhauser@maths.ox.ac.uk \\
       \addr Mathematical Institute \\
       University of Oxford\\
       Oxford, OX2 6GG, United Kingdom}

\editor{Kenji Fukumizu}

\maketitle
\begin{abstract}%
  The sequence of moments of a vector-valued random variable can characterize its law. 
  We study the analogous problem for path-valued random variables, that is stochastic processes, by using so-called robust signature moments.
  This allows us to derive a metric of maximum mean discrepancy type for laws of stochastic processes and study the topology it induces on the space of laws of stochastic processes.
  This metric can be kernelized using the signature kernel which allows to efficiently compute it.
  As an application, we provide a non-parametric two-sample hypothesis test for laws of stochastic processes.
\end{abstract}
\linepenalty=10000

\section{Introduction}
Many inference tasks about data in a topological space $\calX$ involves
\begin{enumerate}[label=(\Alph*)]
\item\label{learn function} making inference about a function $f \in \calF$ from a class $\calF\subset \RR^\calX$,
\item\label{learn measure} making inference about the probability measure $\mu$ from which the data was sampled.
\end{enumerate}
A common approach is to use a ``feature map'', $\Phi:\calX \rightarrow E $, that maps $\calX$ into a linear space $E$ which is rich enough so that both tasks, \ref{learn function} and \ref{learn measure}, can be addressed by linear methods in the ``feature space'' $E$. 
\begin{example}[vector-valued data]\label{ex:vector data}
  If $\calX \subset \RR^d$ is compact, the map
  \begin{align}
    \Phi: \calX \rightarrow \TSeries{(\RR^d)}, \quad x \mapsto \Big(\frac{x^{\otimes m}}{m!}\Big)_{m \ge 0}
  \end{align}
has the property that   
\begin{enumerate}[label=(\alph*)]
\item\label{itm: vector universal} any $f \in \calF\equiv C(\calX,\RR^d)$ can be arbitrary well-approximated by a
  linear functional of $\Phi$,
\item \label{itm: vector characteristic} 
  the map $\mu \mapsto \EE[\Phi(x)] \in \TSeries{(\RR^d)}$, which maps a probability measure $\mu$ on $\calX$ to a sequence of tensors, is injective.
\end{enumerate}
Here,~\ref{itm: vector universal} follows directly from the Stone--Weierstrass theorem which guarantees that linear combinations of monomials approximate continuous functions arbitrary well;~\ref{itm: vector characteristic} is the well-known fact that the moments of a bounded random variable characterize its law.    
\end{example}
This article studies these two questions,~\ref{learn function} and \ref{learn measure}, for the case that $\calX$ is a space of paths, \[\calX \subset \bigcup_{t} C([0,t],\RR^d)\;.\]
This includes the case of sequence-valued data since any sequence of vectors $x(t_1),\ldots,x(t_L) \in \RR^d$ can be identified as a path $x \in C([0,t_L],\RR^d)$ by piecewise linear interpolation between the points $x(t_1),\ldots,x(t_L)$. 
However, first understanding the general continuous time case guarantees that the resulting approach is robust and well-defined for sequence-valued data in the high-frequency limit when the time-discretization gets finer; see also Remark~\ref{rem:sequences}.
Moreover, this covers important classes of processes such as stochastic differential equations or semimartingales and allows us to use tools from stochastic analysis.
By using embeddings into linear spaces, such as through reproducing kernels, we are further able to handle paths evolving in general topological spaces.

\paragraph{Path-valued data.}
A classical theme from stochastic analysis is that for a path $x=(x(t))_{t\in[0,T]}$, the ``signature'' map 
\begin{equation} \label{eq:sigmap}
x \mapsto  \left(\int dx^{\otimes m}\right)_{m \ge 0} \in \TSeries{(\RR^d)}
\end{equation}
is an injection up to tree-like equivalence\footnote{\label{fn:tree-like}Tree-like equivalence is a very useful equivalence relation on pathspace; e.g.~it identifies paths up to a time-change.
  From an analytic point of view, tree-like equivalence is the analogous notion of Lebesgue almost sure equivalence of sets in $\RR^d$ on pathspace.
  We give a precise definition in Section~\ref{subsec:trees}.}
and the tensor $\int dx^{\otimes m} \in (\RR^d)^{\otimes m}$ behaves in a precise algebraic sense analogous to a monomial $\frac{y^{\otimes m}}{m!} \in (\RR^d)^{\otimes m}$ of degree $m$ of a vector $y \in \RR^d$; here we use the shorthand notation $\int dx^{\otimes m}$ for the iterated integral
\begin{align}
  \int dx^{\otimes m}& \coloneqq \int_{0 < t_1 < \cdots < t_m < T}dx({t_1})\otimes \cdots \otimes dx({t_m})\\
                     &\coloneqq \int_{0 < t_1 < \cdots < t_m < T}\dot x({t_1})\otimes \cdots \otimes \dot x({t_m}) dt_1 \cdots dt_m\in (\RR^d)^{\otimes m},\label{eq:signature derivatives} 
\end{align}
where $\dot x(t)\coloneqq (\frac{dx^1(t)}{dt},\ldots,\frac{dx^d(t)}{dt}) \in \RR^d$.
The iterated integrals are defined through either classical Riemann-Stieltjes integrals, or as stochastic or rough path integrals if $x$ is not smooth.
It is then a natural question whether the expected value of these iterated integrals characterizes the law of random paths, analogous to Example~\ref{ex:vector data}.  
\begin{example}\label{ex:BM}
If $Y\sim N(0,1)$, the standard normal on $\RR$, then
\[
\Big(\EE\Big[\frac{1}{m!}Y^{\otimes m}
\Big]\Big)_{m \geq 0}
=
\exp \Big(\frac{1}{2} e_1\otimes e_1\Big) \in \TSeries{\RR},
\]
where $e_1$ is the unit basis vector of $\RR$ and $\exp$ denotes the natural generalisation of the exponential map, see Appendix~\ref{sec:tensor algebra}.
The classical Carleman's condition implies that this moment sequence uniquely determines the law of $Y$.
A Brownian motion $X=(X_t)_{t \in [0,T]}$ in $\RR^d$ is arguably the natural infinite-dimensional analogue of the normal distribution; indeed one has the analogous formula~\citep{Fawcettthesis}
\begin{equation}\label{eq:BM_esig}
\Big( \EE\Big[\int dX^{\otimes m}\Big]\Big)_{m \geq 0} = \exp\Big(\frac{T}{2}\sum_{i=1}^d e_i\otimes e_i\Big)
\in \TSeries{(\RR^d)}\;,
\end{equation}
where $\{e_1,\ldots,e_d\}$ is an orthonormal basis of $\RR^d$ and the stochastic integrals are in the Stratonovich sense.
It was shown in~\citep{ChevyrevLyons16} that~\eqref{eq:BM_esig} uniquely determines the law of Brownian motion (up to tree-like equivalence).
\end{example}

Example~\ref{ex:BM} motivates the signature of a path (resp.~expected signature of a stochastic process) as a generalisation of polynomials (resp.~moments) to path-valued data.
This viewpoint is well-known in stochastic analysis but has become increasingly popular in machine learning; we refer to~\citep{CK16} for a short overview and to~\citep{lyons-04} for a concise mathematical introduction.
We emphasize that such sequences of iterated integrals go at least back to the 20th century and certain aspects of it are studied by various communities in algebra, topology, and control theory (for example, the terms time-ordered exponential, chronological calculus, and Chen-Fliess series all refer to essentially the same object that we call signature; see~\citep[Section 1.4(iii)]{OK16} for such references). 
\begin{remark}
It is natural to ask that, if Example~\ref{ex:BM} generalizes the moment sequence of Example~\ref{ex:vector data} from vector-valued to path-valued random variables (i.e.~the moment generating function), then what is the natural generalization of the Fourier transform (i.e.~the characteristic function)?
Indeed, \citet{ChevyrevLyons16} introduce and study such a ``Fourier transform'', but the resulting objects are quite abstract involving unitary representations of the signature group; in particular it is not clear how to turn this theory into computable features or kernels.
Nevertheless, we believe this to be an interesting research question since already for vector-valued data, both approaches---moment sequences and Fourier transforms---are useful, although arguably moment sequences are more popular for inference tasks.
\end{remark}
\begin{remark}
There are other ways to characterize the law of a stochastic process.
Arguably, the most popular approach is to use finite-dimensional time marginals (via Kolmogorov's extension theorem).
However, this heavily depends on time-parametrization, whereas signature moments are parametrization invariant.
Finally, other kernels for times series exist such as~\citep{chan2005probabilistic, berndt1994using, cuturi2011autoregressive,moreno2003kullback} as do methods coming from dynamic-warping; we refer to \citep[Section 1.4 and Section 5]{OK16} for an attempt of an overview of this large literature.  
What distinguishes the signature kernel is that (i) it applies to general continuous time processes with discrete time arising as a special case---this ensures that it behaves well under high-frequency sampling (see Remark~\ref{rem:sequences})---and (ii) it allows to use strong theoretical results from stochastic analysis which is essential for the proofs of our main theorems.
\end{remark}
\paragraph{What's Wrong With (Signature) Moments.}
Already the feature map from Example~\ref{ex:vector data} runs into trouble if one replaces the compact subset $\calX$ of $\RR^d$ with a non-compact subset of $\RR^d$: in this case, neither~\ref{learn function} or \ref{learn measure} hold in general; a well-known example is that a log-normal random variable is not characterized by its moments. 
Unfortunately this applies, a fortiori, to path-valued data and signature moments: even the laws of ``trivial'' stochastic processes that have straight lines as sample trajectories are in general not characterized by their signature moments; see Appendix~\ref{ex:example nonchar} for a concrete example.   

Even worse, for vector-valued random variables a compact state space can be a reasonable assumption, but this is no longer the case for path-valued random variables since the space of paths is not even locally compact.
A more reasonable alternative is to derive moment decay conditions that ensure that the sequence of moments still characterizes the law despite a non-compact support analogous to how a normally distributed (vector-valued) random variable is characterized by its moments.
This was done in~\citep{ChevyrevLyons16} for signature moments, see Example \ref{ex:BM}.
However, these signature moment decay conditions firstly lead to strong assumptions that are hard to verify, in particular making them unsuitable for a non-parametric machine learning approach; secondly, such decay conditions exclude many important processes such as geometric Brownian motion.

Another related drawback of moments (classical and signature moments) is that of statistical robustness: already for $\RR^d$-valued random variables, the mean $\mu \mapsto \EE_{Y \sim \mu}[Y^{\otimes 1}]=\int y \mu(dy)$ is not robust to outliers, e.g.~for mixtures of the form $\mu= (1-\epsilon) \nu + \epsilon\delta_y$ the Dirac measure at $y$ has too much influence, and the median would be a better ``order $m=1$'' statistic; the same happens for higher moments $m>1$ and such considerations motivated the development of ``robust statistics'', see~\citep{Hampel_robust, huber}.

\subsection{Contribution}
Our main contribution is to derive a ``normalization'' $\lambda(x)$ so that the \emph{robust signature} 
\begin{align}\label{eq:robust sig}
 x\mapsto \Phi(x)=\Big(\lambda(x)^m\int dx^{\otimes m} \Big)_{m\geq0} 
\end{align}
\begin{enumerate}[label=(\alph*)]
\item\label{itm:goodfeature} is a robust, characteristic and universal feature map on pathspace; more precisely 
  \begin{enumerate}[label=(\roman*)]
  \item 
    \emph{universal} means that any continuous bounded function $f(x)$ of paths can be approximated by linear functionals of $\Phi(x)$,
  \item \emph{characteristic} means that the map $\mu \mapsto \EE_{X \sim \mu}[\Phi(X)]$, which maps the law of a stochastic process $X$ to the expected value of $\Phi(X)$, is injective. 
    \item \emph{robust} means that the influence function (in the sense of Huber) of $\Phi$ is finite.
  \end{enumerate}

\item\label{itm:goodk} gives rise to a characteristic and universal kernel $\kernel(x,y)=\langle \Phi(x),\Phi(y) \rangle$ for paths, 
\item\label{itm:metric} the associated maximum mean discrepancy (MMD) is a metric for laws of stochastic processes that yields a topology comparable to the topology of weak convergence under natural assumptions.
\end{enumerate}
The need for such kernel mean embeddings of laws of stochastic processes has been pointed out in the recent survey \citep[Section 6.2]{muandet2017kernel}: 
{\begin{quote}\itshape
``Highly structured data can be taken into account using a kernel, but it is not clear how to employ kernel mean embedding for stochastic processes-i.e., distributions over infinite dimensional objects (...).''  
\end{quote}}
Our robust signature kernel and MMD provides an answer to this question; additionally, it provides the option to ignore the time-parameterization of stochastic processes, which can make the learning much more efficient by quotienting out invariance under a very large space of transformations.
Some of our results remain interesting without the normalization $\lambda$, see Remark \ref{rem:notrobust}.
Below we give a brief outline of our approach.
\paragraph{Robust Signatures.}
A recent observation is that universality and characteristicness are in duality and therefore often equivalent, see~\citep{SGS18}.\footnote{The subtlety is that probability measures form a convex space but not a linear space, so one has to work with general distributions rather than probability measures for characteristicness.}
Hence, it is sufficient to show universality, but this is more involved than in the finite-dimensional case due to the non-local compactness of pathspace. 
To address this, we use a convenient generalization of the Stone--Weierstrass theorem due to~\citet{giles1971generalization}
together with a carefully constructed normalization map $\lambda$.
We refer to the map $\Phi$, as defined in \eqref{eq:robust sig}, as \emph{the robust signature} for two reasons: firstly, and most importantly, it is robust in the sense that no assumptions on the probability measure $\mu$ are needed to ensure characteristicness (except that the iterated integrals are well-defined).
Secondly, the resulting signature moments $\EE_{X \sim \mu}[\Phi(X)]$ that characterize the law $\mu$ of $X$ are robust in the sense of statistical robustness to outliers mentioned above, \citep{Hampel_robust}.

\paragraph{A MMD for Laws of Stochastic Processes.}
A natural distance between probability measures $\mu$,$\nu$ on a topological space $\calX$ is the  ``maximum mean discrepancy'' (MMD)
\begin{align}\label{eq:MMD}
  d\left( \mu,\nu \right)=\sup_{f} \norm{\EE_{X\sim \mu}[f(X)]-\EE_{Y\sim \nu}[f(Y)]}\;,
\end{align}
where the $\sup$ is taken over a sufficiently large set of functions from $\calX$ to $\RR$.
To approximate~\eqref{eq:MMD} from finite samples of the laws $\mu$ (resp.~$\nu$), naive approaches are troublesome because of the supremum over a large space of functions.
To address this, we follow~\citep{berlinet2011reproducing,gretton2012kernel,sriperumbudur2010hilbert} and kernelize the robust signature feature map $\Phi$. 
That is, we define the normalized signature kernel 
$
\kernel(x,y) \coloneqq\langle \Phi(x),\Phi(y) \rangle
$
and take the $\sup$ in~\eqref{eq:MMD} over functions in the unit ball of the reproducing kernel Hilbert space (henceforth RKHS) of $\kernel$.
This gives the identity 
\begin{equation}\label{eq:KMMD}
d^2_{\kernel}(\mu,\nu)=\EE[\kernel(X,X^\prime)]-2\EE[\kernel(X,Y)]+\EE[\kernel(Y,Y^\prime)]\;,
\end{equation}
where $X^\prime$ and $Y^\prime$ denote independent copies of $X\sim \mu$ and $Y\sim \nu$.
Efficient recursive algorithms for inner products of signatures have been developed, see~\citep{OK16}, and we can extend this to the robust signature $\Phi$.
This allows to evaluate $d_{\kernel}(\mu,\nu)$ from finite samples, even if the paths $X$, $Y$ evolve in high-dimensional spaces (large $d$) or general topological spaces. 
Unfortunately, topologies induced by MMDs are typically hard to relate to standard topologies, see~\citep{rachev1991probability}. 
However, one of our main results is that $d_{\kernel}$ is a metric that induces a topology on probability measures strictly weaker than classical weak convergence\footnote{This provides a natural example of a bounded, continuous, characteristic kernel which does not metrize weak convergence on a Polish space.
The existence of such a kernel was raised as an open question in~\citep{SGS18}.
During the writing of this article, this question was additionally answered in~\citep{simongabriel2020metrizing}
for even a locally compact Polish space.
Additionally, the latter corrected an error in~\citep{SGS18} with regards to characterising kernels which metrize weak convergence.} (also called narrow convergence/convergence in law) and give conditions when these topologies are equal.
\paragraph{Application: Hypothesis Testing for Laws of Stochastic Processes.}
We apply our theoretical results to the problem of two-sample hypothesis testing for stochastic processes.
A two-sample test for a stochastic processes $X=(X_t)_{t\in
[0,T]}$, $Y=(Y_t)_{t\in[0,T]}$, tests the null-hypothesis 
\begin{align*}
  H_0:P_X=P_Y\text{ against the alternative } H_1:P_X\neq P_Y,
\end{align*}
where $P_X \coloneqq\PP\circ X^{-1}$ and $P_Y \coloneqq\PP\circ Y^{-1}$ are the laws of $X$ and
$Y$.
Our results about the characteristicness allow us to use the framework of \citet{gretton2009fast, gretton2012kernel} for this problem. 
The data we use is a popular archive of multi-variate time series data that contains heterogeneous datasets from various domains such as speech, handwriting, motion, etc.
We benchmark against classical tests for multi-variate data and MMD tests using classical kernels. 

\subsection{Outline}
Section~\ref{sec:univ. and char.} introduces the equivalence between universality and characteristicness of a feature map $\Phi:\calX \rightarrow E$, as well as the strict topology that becomes useful when $\calX$ is a general non-locally compact space.
Section~\ref{subsec:normalization} focuses on the tensor algebra, $E \coloneqq\FSeries{V}$, as the natural feature space for moments and introduces the concept of a tensor normalization $\Lambda:\FSeriesOne V \rightarrow \FSeriesOne V $. 
As an application, we revisit Example~\ref{ex:vector data} and show how such a normalization turns classical monomials on $\calX=\RR^d$ into a robust feature map for $\RR^d$-valued random variables.
Sections~\ref{sec: ordered moments} and \ref{sec:rough paths} introduce the robust signature map which extends this argument from $\RR^d$-valued random variables to path-valued random variables. 
Section~\ref{sec:kernel} introduces a bounded, continuous, universal and characteristic kernel $\kernel(x,y)=\langle  \Phi(x) , \Phi(y) \rangle$ for paths $x,y$
and discusses the associated MMD on the space of laws of stochastic processes, the topology it induces,
and its discretisations and computability.
Section~\ref{sec:twosample} contains as application a two-sample test between laws of stochastic processes.
The Appendices contain technical details on rough paths, kernels, and the conducted experiments.

\begin{remark}\label{rem:sequences}
\citet{OK16} introduced discretized signature kernels for sequence-valued data and provide efficient algorithms with a focus on supervised classification.
We emphasize that, even for sequence-valued data, the question of characteristicness is non-trivial due to the lack of local compactness.
This is even more the case for path-valued data, which we consider here.
The move from sequence-valued to path-valued data is not only of theoretical interest: although real-world data sets often consist of sequences $(x({t_i}))_{i=1,\ldots, l}$ rather than paths $(x(t))_{t \in [0,T]}$, understanding the limiting behaviour as $\max_i|t_{i+1}-t_i|\rightarrow0$ is crucial; e.g.~it guarantees that no constants blow up in the high-frequency limit.
\end{remark}

\begin{remark}\label{rem:notrobust}
  If one ignores the normalization, that is replaces $\lambda(x)$ by $1$ in \eqref{eq:robust sig}, then the signature features and kernel are in general not universal or characteristic; for example, the expected signature fails to characterize even some trivial stochastic processes, see~Appendix~\ref{ex:example nonchar}.
  However, their use might still be justified in certain situations (for example, if additional moment decay assumptions are known) and many of our results remain interesting in this case, e.g.~the (non-robust) signature MMD can still be used for hypothesis testing, albeit under much stronger assumptions.    
  Nevertheless, the computational cost of normalization is negligible and the robust signature (kernel) enjoys beyond universality and characteristicness also B-robustness which is in particular relevant for real-world data for the same reasons why robustness matters for vector-valued data; hence, in general we recommend the use of the robustified signature kernel resp.~MMD. 
\end{remark}
\subsection{Notation}

We collect some commonly used notation in the following table.

\begin{center}

 \renewcommand{\arraystretch}{1.1}
 \begin{longtable}{lll}\toprule
 Symbol & Meaning & Page\\
 \toprule
      \multicolumn{2}{c}{Spaces}
      \\
      \midrule
   $\calX$ & topological space \\
 $E$ & topological vector space (TVS) over $\RR$  \\
 $E^*$ & algebraic dual of $E$, i.e., space of linear functionals $E \to \RR$ \\
   $E'$ & topological dual of $E$, i.e., subspace $E'\subset E^*$ of \emph{continuous} functionals \\
$H$ & Hilbert space over $\RR$\\
$V$ & Banach space over $\RR$\\
$\overline V$ & $\overline V = V\oplus\RR$ & \pageref{overline_V_page_ref}\\
      \midrule
      \multicolumn{2}{c}{Paths and sequences}\\
      \midrule
   $\pathspVar{1}$ & the subset of $\cup_{T>0}C([0,T],V)$ of absolutely continuous paths & \pageref{page ref 1var RRd}\\
$\rho$ & pseudo-metric on $\pathspVar{1}$ & \pageref{eq:rho}\\
$\overline\rho$ & metric on $\pathspVar{1}(\overline V)$ & \pageref{overline_rho_pageref}\\
   $\pathsreduced{1}$ & tree-like equivalence classes of absolutely continuous paths in $V$& \pageref{def:BV tree-like equivalence}\\
   $x \sim_t y$ & tree-like equivalence relation between $x,y \in \cup_{T>0}C([0,T],V)$ & \pageref{page ref treelike} \\
      $\pi$ & partition of $[0,T]$, i.e.,~a collection of points $0 \le t_1<\cdots <t_l \le T$ \\
      $x^\pi$ &  the sequence $x^\pi\coloneqq(x(t_i))_{i=1,\ldots,l}$ given by sampling $x$ along $\pi$ \\
\midrule
\multicolumn{2}{c}{Signatures and normalization}\\
\midrule
   $\int dx^{\otimes m}$ & shorthand for $\int_{0\leq t_1 \leq \cdots \leq t_m \leq T} dx(t_1) \otimes \cdots \otimes dx(t_m)$ & \pageref{eq:mthIterInt} \\
 $\Sig$ & signature map, $\Sig(x)=\left(\int dx^{\otimes m}\right)_{m\geq0}$ & \pageref{eq:ordered moments map BV}  \\
   $\FSeries{V}$& square summable elements in the tensor algebra $\TSeries{V}$&\pageref{page ref FSeries}\\
$\FSeriesOne{V}$& subset of $\FSeries{V}$ with zero-th component equal $1$ &\pageref{page ref FSeriesOne}\\
$\Lambda$ & a tensor normalization & \pageref{page ref tensor norm}\\
 \bottomrule
 \end{longtable}
\end{center}

\section{Learning in Non-locally Compact Spaces}\label{sec:univ. and char.}
As discussed in the Introduction, two classical problems of learning in a topological space $\calX$ are: \ref{learn function} make inference about a function $f\in\calF\subset \RR^\calX$, and~\ref{learn measure} make inference about a probability measure $\mu$ on $\calX$. 
A standard approach is to map $\calX$ into a (typically infinite or high dimensional) linear space $E$ and address the learning problem there by using linear methods.
All vector spaces in the sequel are over $\RR$.

\begin{definition}\label{def:char. and universal features}
  Let $\calX$ be a topological space.
  For a topological vector space (TVS) $E$, we call any map
$
    \Phi:\calX \rightarrow E
$
  a feature map for which $E$ is the feature space.
  We denote by $E'$ and $E^*$ the topological and algebraic duals of $E$ respectively.
\end{definition}

\paragraph{Universal and Characteristic Features.}
To address Point~\ref{learn function}, we require that $E^\prime$ is large enough to approximate the elements of the function class $\calF \subset \RR^\calX$, i.e.~for every $f\in\calF$, there exists $\ell\in E^\prime$ such that 
$
  f(\cdot)\approx \langle \ell,\Phi(\cdot) \rangle
$
as real-valued functions on $\calX$.
To address Point~\ref{learn measure}, we require that the feature map is non-linear enough to distinguish measures $\mu$ on $\calX$, i.e.~the map
\[
  \mu \mapsto \mu(\Phi)\coloneqq\int_\calX \Phi(x)\mu(dx) \in E
\]
is injective. 
Since $E$ can be infinite dimensional, it is more convenient to write this map as
\[
  \mu \mapsto \Big(\ell \mapsto \int_\calX
    \langle\ell,\Phi(x)\rangle\mu(dx)\Big) \in \left( E^{\prime} \right)^*\;.
\]
\label{page ref E' Estar}More generally, we can replace
the integral by any distribution $D \in \calF^\prime$, i.e.~require that 
\[
  D \mapsto (\ell \mapsto D(x\mapsto \langle \ell, \Phi(x)\rangle)
\]
is injective.
The definition of universality and characteristicness of $\Phi$ makes this precise.
\begin{definition}
  Fix $\calX$ and a TVS $\calF\subset \RR^\calX$.
  Consider a feature map
$\Phi:\calX \rightarrow E$
such that $\langle \ell ,\Phi(\cdot)  \rangle\in \calF$ for all $\ell \in E^\prime$.
  We say that $\Phi$ is
  \begin{enumerate}[label=(\alph*)]
  \item universal to $\calF$ if the following map has a dense image:
    \begin{equation}\label{eq:iEmbedding}
    \imath: E^\prime\to \calF, \quad \ell \mapsto \langle \ell ,\Phi(\cdot)  \rangle\;.
    \end{equation}
  \item characteristic to a subset $\calP \subset \calF^\prime$ if the following map is injective:
    \[
      \kappa:\calP \rightarrow (E^{\prime})^*,\quad D\mapsto [\ell
      \mapsto D(\langle \ell ,\Phi(\cdot)  \rangle) ]\;.
    \]
  \end{enumerate}
\end{definition}
Universality directly addresses Point~\ref{learn function};
characteristicness addresses Point~\ref{learn measure} in a much more general sense since the dual $\calF^\prime$ is larger than the set of probability measures on $\calX$.
This generalization allows to work on the linear space $\calF^\prime$ (instead of a convex subset)
and allows for the following simple statement of duality.
\begin{theorem}\label{thm:univIffChar}
Suppose that $\calF$ is a locally convex TVS. A feature map $\Phi$ is universal to $\calF$ iff $\Phi$ is characteristic to $\calF^\prime$.
\end{theorem}
The proof of Theorem~\ref{thm:univIffChar} employs the Hahn--Banach theorem and is identical to that of a corresponding result for kernels~\citep[Thm.~6]{SGS18}.
We will see that this result is useful since universality is often much easier to show than characteristicness.
\begin{example}[Example~\ref{ex:vector data} revisited]\label{Ex:compactMoments}
Here, $\calX\subset \RR^d$ compact, $E \coloneqq \prod_{m\geq0}\left(\RR^d\right)^{\otimes m}$ and $\Phi(x) \coloneqq\Big(\frac{x^{\otimes m}}{m!}\Big)_{m\geq 0}$.
The image $\imath(E')$ is the space of polynomials on $\calX$.
By Stone--Weierstrass, it follows that $\Phi$ is universal to $\calF \coloneqq C(\calX,\RR)$ equipped with the uniform topology.
By Theorem~\ref{thm:univIffChar}, $\Phi$ is characteristic to $\calF'$ which is the space of regular Borel measure on $\calX$.
In particular, this implies the well-known result that the law of a bounded random variable is determined by its moments.
\end{example}
\subsection{The Strict Topology}\label{sec:strict}
To extend Example~\ref{Ex:compactMoments} from a compact subset $\calX$ of $\RR^d$ to a general non-compact data domain $\calX$, we need to  
\begin{enumerate}[label=(\alph*)]
\item\label{itm:mon} find a map $\Phi:\calX\rightarrow E$ that behaves like  ``monomials'' so that linear functionals and consequently linear functionals $\langle \ell,\Phi(x)\rangle$ behave like ``polynomials'' in $x$; put more abstractly, $ \calF_0 \coloneqq \{x \mapsto \langle \ell, \Phi(x) \rangle: \ell \in E'\}$ should form a point-separating algebra, %
\item\label{itm:F} find a TVS $\calF\subset \RR^\calX$ such that $\calF^\prime$ includes the probability measures on $\calX$,
\item\label{itm:SW} have a Stone--Weierstrass-type result that shows that $\calF_0$ is dense in $\calF$.
\end{enumerate}
Point~\ref{itm:mon} must take into account the concrete choice of $\calX$.
Unfortunately, already for the arguably simple case $\calX=\RR^d$ we know that the standard monomial map from Examples~\ref{ex:vector data} and~\ref{Ex:compactMoments} does not work since there exist unbounded random variables that have different laws but the same moment sequence, see Appendix~\ref{ex:example nonchar}.
Further, for non-locally compact spaces $\calX$ (like the space of paths) it is often easy to find a function class $\calF$ such that exactly one of~\ref{itm:F} and~\ref{itm:SW} holds, but not both; typically the dual space $\calF'$ is either too small, or Stone--Weierstrass-type results are not known or involve conditions that are hard to verify.

Below we recall a result of \citet{giles1971generalization} which provides a surprisingly elegant and general roadmap on what is needed to addresses all points~\ref{itm:mon},~\ref{itm:F}, and~\ref{itm:SW}. 
\begin{definition}
Let $\calX$ be a topological space.
We say that a function $\psi:\calX\rightarrow \RR$ vanishes at infinity if for all $\epsilon>0$ there exists a compact set $K\subset\calX$ such that $\sup_{x\in\calX\setminus K}|\psi(x)|<\epsilon$.
Denote with $B_0(\calX,\RR)$ the set of functions that vanish at infinity. 
The strict topology
on $C_b(\calX,\RR)$ is the topology generated by the seminorms
\[
p_\psi(f)=\sup_{x\in\calX}|f(x)\psi(x)|\;,\quad \psi\in B_0(\calX,\RR)\;.
\]
\end{definition}
\begin{theorem}[\citet{giles1971generalization}]\label{thm:strict top}
Let $\calX$ be a metrizable topological space.
  \begin{enumerate}[label=(\roman*)]
  \item\label{itm:comparison} The strict topology on $C_b(\calX,\RR)$ is weaker than the uniform topology and stronger than the topology of uniform convergence on compact sets.
  \item\label{itm:Stone-Weierstrass strict top} If $\calF_0$ is a subalgebra of $C_b(\calX,\RR)$ such that
    \begin{enumerate}
    \item $\forall x\neq y\in \calX$ there exists $f\in\calF_0$ such that $f(x)\neq
      f(y)$, and
      \item $\forall x \in \calX$ there exists $f\in \calF_0$ such that
        $f(x)\neq 0$,
    \end{enumerate} 
    then $\calF_0$ is dense in $C_b(\calX,\RR)$ under the strict topology.
    \item\label{itm:dual strict top} The topological dual of $C_b(\calX,\RR)$ equipped with the strict topology is
the space of finite regular Borel measures on $\calX$.
  \end{enumerate}
 \end{theorem}
  \begin{proof}
Point~\ref{itm:comparison} follows from the definition.
Point~\ref{itm:Stone-Weierstrass strict top} is~\citep[Thm.~3.1]{giles1971generalization}; in fact, the result holds more generally for
any topological space and one only needs that point-separating and
    non-vanishing functions are in the closure of $\calF_0$.
   Point~\ref{itm:dual strict top} is \citep[Thm.~4.6]{giles1971generalization}.
  \end{proof}
  The specific construction of $\Phi$ must depend on the domain $\calX$.
  However, often one has a natural candidate for monomials on $\calX$.
  For example, for $\calX=\RR^d$ one can use the standard monomials as used in Examples~\ref{ex:vector data} and~\ref{Ex:compactMoments};
  for $\calX$ the space of paths, stochastic analysis suggests that the iterated integrals $\int dx^{\otimes m}$ are the ``natural monomials'' of degree $m$. 
  It remains to have a method that generically turns such monomials into ``bounded monomials'' to which Theorem~\ref{thm:strict top} applies. 
  In Section \ref{subsec:normalization} we introduce a generic normalization that accomplishes this.
  In Section~\ref{sec: ordered moments} we apply this construction to the space of paths to derive our robust signature features.
  \section{Robustification by Normalization}\label{subsec:normalization}
  Already the standard monomials in $\RR^d$, Example~\ref{ex:vector data}, are unbounded on non-compact domains $\calX$, hence linear functionals of them (polynomials) are not bounded and thus Theorem~\ref{thm:strict top} does not apply. 
  On the other hand, polynomials have many nice properties: they are closed under multiplication and they separate points.  
  In this section we show that there exists a generic ``normalization'' that makes monomials bounded while preserving their algebra and point-separating property.
  Moreover, this normalization also turns them into robust statistics in the sense of \citet{huber}. 
  Even for the case $\calX=\RR^d$ this leads to a novel moment sequence that characterize the law of any $\RR^d$-valued random variable, Section~\ref{sec:moments}. 
  More importantly for our purposes, this method also applies to signature features on paths. 
\subsection{The Tensor Algebra}\label{sec:normFunc}
Throughout the rest of this article let us denote with $V$ a Banach space.
Recall that $V^{\otimes m}$ denotes the space of tensors of degree $m$.
A \emph{sequences of tensors of increasing degree} is an element of the space $\TSeries{V}$.
We will denote elements $\bft\in \TSeries{V}$
by
\begin{align*}
  \bft=(\bft^m)_{m\ge 0} \in \TSeries{V}
\end{align*}
where $\bft^m\in V^{\otimes m}$.
The space $\TSeries{V}$ becomes a linear space by defining 
\begin{align*}
  \bfs +\bft & \coloneqq \left( \bfs^0+\bft^0,\bfs^1+\bft^1,\ldots \right)\text{ for }\bfs=(\bfs^m)_{m \geq 0},\,\bft=(\bft^m)_{m\geq 0}\in \TSeries{V}.
\end{align*}
This space is the natural state space for monomials; for general background on tensors we refer to Appendix \ref{sec:tensor algebra}. 
Of particular relevance is the subspace $\FSeries{V} \subset \TSeries{V}$ of elements that have a finite norm. 
\begin{definition}
  Let $V$ be a Banach space.\footnote{We implicitly assume that $V^{\otimes m}$ for all $m \geq 2$ are Banach spaces formed by the completion with respect to some admissible system of tensor norms $\|\cdot\|_{V^{\otimes m}}$, see~\citep[Sec.~3.1]{lyons-qian-02}.}
  We denote by $\FSeries{V}$\label{page ref FSeries} the Banach space
  \[
  \FSeries{V} \coloneqq \Big\{\bft \in \TSeries{V}\, :\, \|\bft\| \coloneqq \Big(\sum_{m\geq0} \|\bft^m\|^2_{V^{\otimes m}}\Big)^{1/2} < \infty \Big\}\;.
  \]
  Define further the subset\label{page ref FSeriesOne} $\FSeriesOne{V} \coloneqq \{\bft \in \FSeries{V}\, :\, \bft^0=1 \}$.
\end{definition}
\begin{remark}
If $V=H$ is a Hilbert space,
we will always equip $H^{\otimes m}$ with the canonical inner product (which induces an admissible system of norms)
given on elementary tensors by
\[
\langle x_1\otimes\ldots\otimes x_m , y_1\otimes\ldots \otimes y_m \rangle_{H^{\otimes m}} = \prod_{j=1}^m \langle x_j,y_j\rangle_H\;,
\]
and extended by linearity to $H^{\otimes m}$.
Consequently, $\FSeries{H}$ becomes a Hilbert space with inner product $\langle \bfs, \bft\rangle_{\FSeries{H}} = \sum_{m \geq 0} \langle \bfs^m,\bft^m\rangle_{H^{\otimes m}}$.
\end{remark}

\subsection{Normalization}
Naive approaches to make monomials bounded such as dividing them by their norm destroys that they form an algebra and that they are point-separating. 
The idea of a tensor normalization, is that these ``nice'' properties of monomials can be preserved while making them bounded.
\begin{definition}\label{def:norm}
A \emph{tensor normalization}\label{page ref tensor norm} is a continuous injective map of the form
  \begin{align*}
 \Normal : \FSeriesOne{V} \rightarrow \left\{\bft \in \FSeriesOne{V}:\, \|\bft\|\leq R\right\}\;,
 \qquad
  \Normal : \bft \mapsto \delta_{\lambda(\bft)} \bft \;,
  \end{align*}
 where $R > 0$ is a constant, $\lambda : \FSeriesOne{V} \to (0,\infty)$ is a function, and, for $c\in\RR$, $\delta_c$ is the dilation map $
\delta_c(\bft) = (\bft^0, c\bft^1, c^2 \bft^2, \ldots )$.
\end{definition}

The existence of tensor normalizations is not trivial.
In the rest of this section, we give a general method to construct such maps and determine their regularity properties.
In Sections~\ref{subsec:comp_sig_kernel} and~\ref{subsec:implementation} we show that tensor normalization can be carried out in a kernel learning setting with minor computational cost.

\begin{lemma}\label{lem:dilationBounds}
Suppose that $\lambda \geq 0$ and $\bft \in \FSeries{V}$ such that $\delta_\lambda \bft \in \FSeries{V}$.
Then $\|\delta_\lambda \bft - \bft\|^2 \leq |\|\delta_\lambda \bft\|^2 - \|\bft\|^2|$.
\end{lemma}

\begin{proof}
We may suppose without loss of generality that $\lambda \geq 1$ (otherwise we swap $\bft$ with $\delta_\lambda \bft$ and $\lambda$ with $\lambda^{-1}$).
Taking derivatives in $\bar\lambda$ at some $\bar\lambda \in [1,\lambda)$,
\[
\frac{d}{d\bar\lambda} \|\delta_{\bar\lambda} \bft - \bft\|^2
= \sum_{m=1}^\infty 2m(\bar\lambda^{2m-1} - \bar\lambda^{m-1})\|\bft^m\|^2_{V^{\otimes m}} \leq \sum_{m=1}^\infty 2m\bar\lambda^{2m-1}\|\bft^m\|^2_{V^{\otimes m}} = \frac{d}{d\bar\lambda} \|\delta_{\bar\lambda} \bft\|^2\;,
\]
where all sums are convergent by the assumption that $\sum_{m=1}^\infty \lambda^{2m}\|\bft^m\|^2< \infty$.
Since both derivatives are positive and the first is bounded above by the second, the conclusion follows.
\end{proof}

\begin{proposition}\label{prop:LipschitzBoundPhi}
Let $\psi : [1,\infty) \to [1,\infty)$ with $\psi(1)=1$.
For $\bft \in \FSeriesOne{V}$, let $\lambda(\bft) \geq 0$ denote the unique non-negative number such that $\|\delta_{\lambda(\bft)} \bft\|^2 = \psi(\|\bft\|)$.
Define $\Normal : \FSeriesOne{V} \rightarrow \FSeriesOne{V}$, $\Normal(\bft) = \delta_{\lambda(\bft)} \bft$.
Denote further $\|\psi\|_\infty = \sup_{x \geq 1} \psi(x)$.
\begin{enumerate}[label=(\roman*)]
\item \label{point:maps into ball} It holds that $\Lambda$ takes values in the set $\{\bft \in \FSeriesOne{V} : \|\bft\| \leq \sqrt{\|\psi\|_\infty}\}$.
\item \label{point:injection} If $\psi$ is injective, then so is $\Lambda$.
\item \label{point:Lipschitz} Suppose that $\sup_{x \geq 1}\psi(x)/x^2 \leq 1$, $\|\psi\|_\infty < \infty$, and that $\psi$ is $K$-Lipschitz for some $K > 0$.
Then
\[
\|\Lambda(\bfs) - \Lambda(\bft)\|
\leq \big(1+K^{1/2}+\sqrt 2\|\psi\|_\infty^{1/2}\big)\big(\|\bfs-\bft\|^{1/2}\vee\|\bfs-\bft\| \big)\;.
\]
\end{enumerate}
\end{proposition}

\begin{proof}
\ref{point:maps into ball} is clear by construction.
For~\ref{point:injection} suppose that $\psi$ is injective and that $\Lambda(\bft) = \Lambda(\bfs)$.
If $\bft = \bfone \coloneqq (1,0,0\ldots)$, then evidently $\bfs=\bft$.
We thus suppose $\bft\neq \bfone$.
By definition of $\Lambda$, it holds that $\bfs = \delta_\lambda\bft$ for some $\lambda \geq 0$.
On the other hand, we also have $\psi(\|\bft\|)=\psi(\|\bfs\|)$.
Since $\psi$ is injective and $\lambda \mapsto \|\delta_\lambda \bft\|$ is strictly increasing by the assumption that $\bft = \bfone$, it follows that $\bfs = \bft$, which proves~\ref{point:injection}.
For~\ref{point:Lipschitz}, note that $\psi(x) \leq x^2$ implies $\lambda(\bft) \leq 1$, and thus $\delta_{\lambda(\bft)}$ is $1$-Lipschitz on $\FSeriesOne{V}$.
Hence
\begin{align*}
\|\delta_{\lambda(\bfs)}\bfs - \delta_{\lambda(\bft)}\bft\| &\leq \|\delta_{\lambda(\bft)}\bfs - \delta_{\lambda(\bft)}\bft\| + \|\delta_{\lambda(\bfs)}\bfs - \delta_{\lambda(\bft)}\bfs\|
\\
&\leq \|\bfs-\bft\| + |\|\delta_{\lambda(\bfs)}\bfs\|^2 - \|\delta_{\lambda(\bft)}\bfs\|^2|^{1/2}
\\
&\leq \|\bfs-\bft\| + |\psi(\|\bfs\|) - \psi(\|\bft\|)|^{1/2} + |\|\delta_{\lambda(\bft)}\bft\|^2 - \|\delta_{\lambda(\bft)}\bfs\|^2|^{1/2}
\\
&\leq \|\bfs-\bft\| + K^{1/2}\|\bfs-\bft\|^{1/2} + \sqrt 2\|\psi\|_{\infty}^{1/2}\|\bfs-\bft\|^{1/2}\;,
\end{align*}
where we used in the second line Lemma~\ref{lem:dilationBounds} and in the fourth line that $\psi$ is $K$-Lipschitz and without loss of generality that $\lambda(\bft) \leq \lambda(\bfs)$ (hence $\|\delta_{\lambda(\bft)}\bfs\| \leq \|\delta_{\lambda(\bfs)}\bfs\| \leq \|\psi\|_\infty$).
\end{proof}

\begin{corollary}\label{cor:tensor_normal}
Let $\psi:[1,\infty) \to [1,\infty)$ be injective satisfying $\psi(1)=1$ and the conditions of Proposition~\ref{prop:LipschitzBoundPhi}\ref{point:Lipschitz}.
Then $\Lambda$ constructed in Proposition~\ref{prop:LipschitzBoundPhi} is a tensor normalization.
\end{corollary}

\begin{example}\label{ex:psi}
The tensor normalization $\Lambda$ we use in our experiments is built as in Corollary~\ref{cor:tensor_normal} from
the function $\psi:[1,\infty)\to[1,\infty)$ defined by
\[
\psi(\sqrt x) = 
\begin{cases}
x &\text{ if } x \leq C\;,
\\
C + C^{1+a}(C^{-a}-x^{-a})/a &\text{ if } x > C\;,
\end{cases}
\]
where $a>0,C\geq 1$ are constants.
Note that $\psi$ satisfies all the conditions of Corollary~\ref{cor:tensor_normal}; in particular,
$\psi$ is $2C$-Lipschitz and bounded above by $C(1+\frac1a)$.
Figure~\ref{fig:psi} shows the graph of $\psi(\sqrt x)$ with $a=1$ and $C=4$; note that
$\Lambda(\bft)=\bft$
if $\|\bft\| \leq 4$.
\begin{figure}[h]
\centering
\begin{tikzpicture}
     \begin{axis}[
       name =Wn,
       height       = 2.0in,
       width       = 4.0in,
       xmax         = 20,
       ymax         = 10,
       xtick        = {0.0, 2, 4, 6,8,10,12,14,16,18,20},
       xticklabels  = {0, $2$, $4$, $6$, $8$, $10$,$12$,$14$,$16$,$18$,$20$},
       axis x line  = bottom,
       ytick        = {0.0, 2, 4, 6,8,10},
       yticklabels  = {0, $2$, $4$, $6$, $8$,$10$},
       axis y line  = left,
       line cap=round
     ]
         \addplot[gray,dotted] coordinates { (0,10) (20,10) };
         \addplot[gray,dotted] coordinates { (20,0) (20,10) };
         \addplot[gray,dotted] coordinates { (0,0) (10,10) };
         \addplot[gray,dotted] coordinates { (4,0) (4,10) };
         \addplot[gray,dotted] coordinates { (0,8) (20,8) };
         \addplot[red,very thick,domain=1:4, samples=2]{ x };
         \addplot[red,very thick,domain=4:20, samples=40]{ 4+16*(0.25-x^(-1)) };
     \end{axis}
     \node[anchor=north] at (Wn.south) {$\begin{matrix} \\ x \end{matrix}$};
         \node[anchor=east] at (-0.5,2) {$\begin{matrix} \\ \psi(\sqrt x) \end{matrix}$};
     \end{tikzpicture}
\vspace{-2ex}
\caption{The graph of $\psi(\sqrt x)$ from Example~\ref{ex:psi} with $a=1$ and $C=4$.}
\label{fig:psi}
\end{figure}
\end{example}

\subsection{B-robust Statistics}
In addition to allowing us to applying Theorem~\ref{thm:strict top},
tensor normalizations guarantee that the statistics which are linear functionals $\langle \ell, \Phi(X) \rangle$ of a normalized feature map $\Phi=\Lambda \circ \varphi$ are robust to outliers in a precise statistical sense.   
\begin{definition}[\citet{Hampel_robust}]
Let $\calP(\calX)$ denote the space of Borel probability measures on $\calX$.
For a map $T: \calP(\calX) \to \RR$ define 
the \emph{influence function} of $T$ at $\mu\in\calP(\calX)$ by
 \[
   \IF(x;T,\mu) \coloneqq \lim_{\epsilon \searrow 0} \frac{T(\epsilon \delta_x + (1-\epsilon) \mu) - T(\mu)}{\epsilon}
\]
for $x\in\calX$ for which this limit exists. $T$ is called \emph{$B$-robust} at $\mu$ if
$\sup_{x\in\calX} |\IF(x;T,\mu)|<\infty$.
\end{definition}

\begin{example}
Suppose $f:\calX\to \RR$ is measurable
and $T(\mu) \coloneqq \mu(f) = \int_{\calX} f(x)\mu(dx)$
if $f$ is $\mu$ integrable (and arbitrary otherwise).
Then, provided $f$ is $\mu$-integrable, $\IF(x;T,\mu) = f(x)-\mu(f)$,
and therefore $T$ is $B$-robust at $\mu$ if and only if $f$ is bounded.
\end{example}

The previous example shows that if $\calX\subset \RR$ and $T:\calP(\calX) \to \RR$ is the $m$-th moment $T(\mu) = \EE_{X\sim\mu}[X^m]$ for some $m \geq 1$,
then $T$ is $B$-robust (at any $\mu$ with finite $m$-th moment) if and only if $\calX$ is bounded.
The same example, however, shows the following.

\begin{proposition}
Consider $\Lambda$ a tensor normalization, $\varphi:\calX \to \FSeriesOne{V}$ measurable, and $\ell\in \FSeries{V}'$. Define $\Phi=\Lambda\circ\varphi$ and $f(x) = \scal{\ell, \Phi(x)}$.
  Then $\mu\mapsto \mu(f)$ is $B$-robust at any $\mu$.
\end{proposition}

The significance of this result is that linear functionals of the feature maps `$\Phi$' appearing in our results below (Proposition~\ref{prop:noncompactMoments}, Theorems~\ref{thm:universal map BV paths} and~\ref{thm:featureMapRPs}, etc.)
are automatically $B$-robust statistics.
The focus of this article is to characterize probability measures on pathspace, but the above observation that $B$-robustness is intimately linked to moment normalization and characteristicness is already interesting for $\RR^d$-valued random variables.

\subsection{Example: Robust Moments for Vector-valued Data}\label{sec:moments}
We now demonstrate the usefulness of the tensor normalization by turning the standard monomial feature map from Example~\ref{Ex:compactMoments} into a robust, characteristic, and universal feature map on the whole of $\RR^d$.
We emphasise that our argument does not rely on local compactness and thus paves the way to apply the same reasoning to pathspace.
\begin{proposition}\label{prop:noncompactMoments}
Consider the map
\begin{equation}\label{monomial map}
\varphi : \RR^d \to \FSeriesOne{\RR^d};, \quad x \mapsto \Big(\frac{x^{\otimes m}}{m!}\Big)_{m=0}^\infty\;.
\end{equation}
Let $\Lambda : \FSeriesOne{\RR^d} \to \FSeriesOne{\RR^d}$ be a tensor normalization. Then the map
  \begin{align*}
    \Phi:  \RR^d\rightarrow \FSeriesOne{\RR^d},
    \quad \Phi=\Normal\circ \varphi
  \end{align*}
  \begin{enumerate}[label=(\roman*)]
  \item \label{point:monomialCont} is a continuous injection from $\RR^d$ into a bounded subset of $\FSeriesOne{\RR^d}$,
  \item \label{point:monomialUniv} is universal to $\calF=C_b(\RR^d,\RR)$ equipped with the strict topology,
  \item \label{point:monomialChar} is characteristic to $\calF^\prime$, the set of finite, signed Borel
    measures on $\RR^d$.
  \end{enumerate}
  Moreover, for any $\ell \in \FSeries{V}'$ the map $\mu \to \langle \ell, \EE_{X \sim \mu} [\Phi(X)] \rangle$ is $B$-robust.
\end{proposition}
\begin{proof}
\ref{point:monomialCont} follows from the definition of the tensor normalization $\Normal$ and the fact that $\varphi : \RR^d \to \FSeriesOne{\RR^d}$ is a continuous injection.
To show~\ref{point:monomialUniv}, we claim that the family of functions
$
\calF_0 = \{x\mapsto\langle \ell, \Phi(x) \rangle, \ell\in (\RR^m)^\prime\,,\; m \geq 0\} \subset C_b(\RR^d,\RR)   
$
  satisfies the conditions of Theorem~\ref{thm:strict top}\ref{itm:Stone-Weierstrass strict top}.
  Indeed, for $m,n \geq 0$ and $\ell \in (\RR^{d})^{\otimes m}, \tilde \ell \in
  (\RR^{d})^{\otimes n}$ we have
$
  \langle\ell,\Phi(x)\rangle \langle\tilde \ell, \Phi(x)\rangle
    =
    \lambda(\varphi(x))^{m+n} \langle\ell,\frac{x^{\otimes m}}{m!}\rangle\langle\tilde\ell,\frac{x^{\otimes n}}{n!}\rangle
 =\binom{m+n}{m}\langle  \ell\otimes\tilde\ell,\Phi(x) \rangle 
$,
  which shows that $\calF_0$ is a subalgebra of $C_b(\RR^d,\RR)$.
The family $\calF_0$ separates points since $\Phi$ is injective.
Finally, note that $\Phi(x) \neq 0$ for all $x \in \RR^d$.
It follows by Theorem~\ref{thm:strict top}\ref{itm:Stone-Weierstrass strict top} that $\Phi$ is universal to $C_b(\RR^d,\RR)$ equipped with the strict topology which shows~\ref{point:monomialUniv};~\ref{point:monomialChar} now follows by Theorem~\ref{thm:strict top}\ref{itm:dual strict top}.
\end{proof}
An informal point of view of Proposition \ref{prop:noncompactMoments} is to recall the difference between the moment generating function and the characteristic function of $\RR^d$-valued random variable $X$.
The former can fail to characterize laws of unbounded random variables, whereas the latter characterizes the law of any random variable.
The characteristic function accomplishes this by using a Fourier transform that spins faster in a circle in the complex plane for large value of $X$.  
The above proposition shows that the same can be accomplished by ``smearing out'' large realizations of $X$ into a ball in $\FSeries{\RR^d}$.   

\section{Robust Features for Smooth Paths}\label{sec: ordered moments}
In this section, we apply the robustification introduced in Section~\ref{subsec:normalization} on pathspace. 
We first focus on the case when $\calX=\pathspVar1$ denotes the subset of $\cup_{T > 0}C([0,T],V)$ of absolutely continuous paths.
The case for more irregular paths is covered in Section~\ref{sec:rough paths}.

\subsection{Iterated Integrals: Monomials of Paths}
We consider first absolutely continuous paths parameterized over any finite interval $[0,T]$ whose starting point is the origin:\label{page ref 1var RRd}
\begin{align*}
  &\pathspVar{1}(V) \coloneqq \bigcup_{T > 0}
\Big\{x \in C([0,T],V): \exists \dot x \in L^1([0,T],V)\;, \; x = \int_0^\cdot \dot x(t) dt \Big\}\;,
\\
  &\|x\|_{1\Var;[0,T]} \coloneqq \int_0^T \|\dot x(t)\|  dt = \sup_{\substack{n \geq 1 \\ 0\leq t_1\leq \cdots \leq t_n \leq T}} \sum_{i=1}^{n-1} \| x(t_{i+1})-x(t_i)\|\;.
\end{align*}
When the space $V$ is clear from the context we just write $\pathspVar{1}$ instead of $\pathspVar{1}(V)$.
When we speak of an element $x\in\pathspVar{1}$, we let $[0,T_x]$ denote the corresponding time interval.  
Due to the arbitrary time interval $[0,T]$, $\pathspVar{1}$ is not a vector space (since two paths on different time intervals can not be added).
To define a topology on $\pathspVar{1}$, consider for $x,y \in \pathspVar{1}$ the quantity\footnote{Readers will notice the similarity with dynamic time warping and Fr{\'e}chet distance.}
\begin{equation}\label{eq:rho}
\rho(x,y) = \inf_{\tau} \|x -y\circ\tau \|_{1\Var;[0,T_x]}\;,
\end{equation}
where the infimum is taken over all increasing bijections $\tau : [0,T_x] \to [0,T_y]$; we call such a $\tau$ a reparameterization.
Note that $\rho$ satisfies all the properties of a metric except that $\rho(x,y) = 0$ does in general not imply that $x=y$ (i.e., $\rho$ is a pseudo-metric).
We equip $\pathspVar{1}$ with the (non-Hausdorff) topology induced by $\rho$.
The $\pathspVar{1}$ regularity of the paths allows us to use Riemann--Stieltjes integration to define (see also~\eqref{eq:signature derivatives}) 
\begin{align}\label{eq:mthIterInt}
  \int_0^t dx^{\otimes m}  \coloneqq\int_{0\leq t_1 \leq \cdots \leq t_m \leq t} dx(t_1) \otimes \cdots \otimes  dx(t_m)\in V^{\otimes m} \quad \text{ for } t\in[0,T_x]\;.
\end{align}
We will often use the shorthand $\int dx^{\otimes m}=\int_0^{T_x} dx^{\otimes m}$.
A guiding principle that plays a prominent role in stochastic analysis is that the so-called \emph{signature map}
\begin{equation}\label{eq:ordered moments map BV}
  \Sig: \pathspVar1 \rightarrow \FSeriesOne{V}\;,\quad
  x \mapsto \Big(1,\int dx, \int dx^{\otimes 2}, \int dx^{\otimes 3}, \ldots\Big)\;,
\end{equation}
is the natural generalization of the monomial feature map from Example~\ref{ex:vector data} to pathspace. 
\begin{example}\label{ex:signature linear path}
  Let $V=\RR^d$. Consider the linear path $x:[0,1] \to \RR^d$, $x(t)=tv$, for a fixed vector $v\in \RR^d$.
  A simple calculation shows $\int dx^{\otimes m}=\frac{v^{\otimes m}}{m!}$.
  Hence, $\Sig$ is a generalization of the monomial feature map on $\RR^d$ from Example~\ref{ex:vector data}. 
\end{example}
An algebraic reason why~\eqref{eq:mthIterInt} can be seen as sequence of monomials on pathspace $\pathspVar{1}$ is that their linear span is closed under multiplication.
The so-called shuffle product makes this precise: considering for simplicity $V=\RR^d$, we have for $\mathbf{i}=(i_1,\ldots,i_m)\in \{1,\ldots,d\}^m,\mathbf{j}=(j_1,\ldots,j_n)\in \{1,\ldots,d\}^n$
  \begin{align}\label{eq:shuffle}
    \int dx^{\otimes m}_{\mathbf{i}}\int dx^{\otimes
    n}_{\mathbf{j}}=\sum_{\mathbf{k}} \int dx^{\otimes (m+n)}_{\mathbf{k}}  \;,
  \end{align}
  where the sum is taken over all $\mathbf{k}=(k_1,\ldots ,k_{m+n})$ that are shuffles\footnote{A shuffle $\mathbf{k}=(k_1,\ldots,k_{m+n})$ of $\mathbf{i}$ and $\mathbf{j}$ is a permutations of $(i_1,\ldots,i_m,j_1,\ldots,j_n)$ subject to the condition that the order of elements in $\mathbf{i}$ and $\mathbf{j}$ is preserved.} of $\mathbf{i}$ and $\mathbf{j}$, see~\citep[Thm.~2.15]{lyons-04}.
The signature map furthermore identifies paths up to tree-like equivalence.
\begin{theorem}[\citet{BGLY16}]\label{thm:propsPhiMap}
  The map $\Sig: \pathspVar{1}\rightarrow \FSeriesOne{V} $ is injective up to tree-like equivalence: for $x,y \in \pathspVar{1}$, $\Sig(x) = \Sig(y)$ if and only if $x\sim_ty$.
\end{theorem}
The tree-like equivalence relation $x \sim_t y$ is a natural generalization of the reparameterization relation on pathspace, see Footnote~\ref{fn:tree-like} and Appendix~\ref{subsec:trees}.
In particular, and perhaps most importantly, two paths $x,y$ are tree-like equivalent if they are reparameterizations, $x(t)=y(\tau(t))$ for all $t$ (or more generally if $\rho(x,y)=0$).
We show in Section~\ref{subsec:timepar BV}, that, if desired, it is easy to remove the tree-like equivalence relation in Theorem~\ref{thm:propsPhiMap}.  
Nevertheless, the identification of tree-like equivalent paths (read: reparameterized paths) can be a powerful dimensionality reduction. 
\subsection{Features for Unparameterized Paths}\label{subsec:univPaths}
The empirical success of dynamic time warping and Fr{\'e}chet distances show that for path-valued data it can beneficial to ignore the time-parameterization. 
Motivated by this we identify paths that are reparameterizations.
Hence, instead of the space of paths $\pathspVar{1}$ we first derive robust features for the space of equivalence classes of paths $\pathsreduced{1}\coloneqq \pathspVar{1} /\sim_t$.
\begin{definition}\label{def:BV tree-like equivalence}
We define the space of unparameterized absolutely continuous paths $\pathsreduced{1}$\label{page ref 1var reduced} as the set of equivalence classes $\pathspVar{1} / \sim_t$.
We equip $\pathsreduced{1}$ with the initial topology\footnote{The quotient map $\pathspVar{1} \to \pathsreduced{1}$, when $\pathspVar{1}$ is equipped with the pseudo-metric $\rho$, is continuous due to the continuity of $\Sig: \pathspVar{1} \to \FSeriesOne{V}$~\citep[Thm.~3.1.3]{lyons-qian-02}.
In particular, the quotient topology on  $\pathsreduced{1}$ is \textit{stronger} than the initial topology we consider.} induced by the signature map embedding $\Sig : \pathspVar{1} \to \FSeriesOne{V}$.
\end{definition}
In the above definition, recall that given a set $\calY$ and a topological space $\calX$, the \emph{initial topology} on $\calY$ induced by a function $f\colon\calY\to\calX$ is the weakest topology on $\calY$ such that $f$ is continuous;
remark that we can view $\Sig$ as a function $\pathsreduced{1} \to \FSeriesOne{V}$ due to Theorem~\ref{thm:propsPhiMap}.

In complete analogy to our toy example of the monomial feature map for $\RR^d$ in Proposition~\ref{prop:noncompactMoments}, we apply a tensor normalization to the signature map $\Sig$ to get a robust, universal, and characteristic feature map.

\begin{theorem}\label{thm:universal map BV paths}
Let $\Lambda : \FSeriesOne{V} \to \FSeriesOne{V}$ be a tensor normalization.
The robust signature
  \begin{equation*}
    \Phi:\pathsreduced{1}\rightarrow \FSeriesOne{V}, \quad \Phi = \Normal \circ \Sig
  \end{equation*} 
  \begin{enumerate}[label=(\roman*)]
  \item \label{point:PhiCont} is a continuous injection from $\pathsreduced{1}$ into a bounded subset of $\FSeriesOne{V}$,
  \item\label{point:PhiUniv} is universal to $\calF \coloneqq C_b(\pathsreduced{1},\RR)$ equipped with the strict
    topology,
  \item\label{point:PhiChar} is characteristic to the space of finite regular Borel measures on $\pathsreduced{1}$.
  \end{enumerate}
  Moreover, $\forall \ell \in \FSeries{V}'$ the map $\mu \to \langle \ell, \EE_{X \sim \mu}[\Phi(X)]\rangle$ is $B$-robust.
\end{theorem}

\begin{proof}
\ref{point:PhiCont} follows from the definition of tensor normalization and the fact that $\Sig: \pathsreduced{1} \to \FSeriesOne{V}$ is a homeomorphism onto its image.
For~\ref{point:PhiUniv}, define $L \coloneqq\bigoplus_{m \geq 0} (V^{\otimes m})'$, which we identify with a dense subspace of $\FSeries{V}'$ via $\ell(\mathbf{t})= \sum_{m \geq 0} \langle\ell^m ,\bft^m\rangle$.
Define further $\calF_0 \coloneqq \{\ell\circ \Phi : \ell \in L\} \subset \calF$.
By point~\ref{point:PhiCont}, $\calF_0$ separates the points $\pathsreduced{1}$.
Furthermore, $\Sig$ takes values in the group-like elements of $\FSeries{V}$~\citep[Cor.~3.9]{CDLL16} which implies $\calF_0$ is closed under multiplication (in the case $V=\RR^d$, this is equivalent to the shuffle product~\eqref{eq:shuffle}). 
It follows that $\calF_0$ satisfies the assumptions of Theorem~\ref{thm:strict top}\ref{itm:Stone-Weierstrass strict top}.
Hence $\calF_0$ is dense in $\calF$ under the strict topology, which proves~\ref{point:PhiUniv}.
Point~\ref{point:PhiChar} in turn follows from Theorem~\ref{thm:univIffChar} and Theorem~\ref{thm:strict top}\ref{itm:dual strict top}.
\end{proof}
\subsection{Features for Parameterized Paths}\label{subsec:timepar BV}
While it can be useful to ignore time-parameterization, there are also many situations where time-parameterization matters, e.g.~financial data.
The simple observation is that we can record the parameterization of a path by adding time as a state-space coordinate. 

Let $\overline V \coloneqq V\oplus\RR$.\label{overline_V_page_ref}
Observe that the map
  \begin{equation*}
    \pathspVar{1}(V) \ni x\mapsto \overline{x} \in \pathspVar{1}(\overline V) \; \text{ with }\;\overline x : [0,T_x] \to \overline V\;,\;\;\overline{x}(t) \coloneqq (x(t),t)
  \end{equation*}
 realises $\pathspVar{1}(V)$ as a subset of $\pathspVar{1}(\overline V)$.
Furthermore, for $x,y \in \pathspVar{1}(V)$, it holds by strict monotonicity of the second component that $\overline x \sim_t \overline y$ if and only if $x=y$.
  Hence, denoting by $\overline\rho(x,y)$\label{overline_rho_pageref} the quantity~\eqref{eq:rho} with $x,y$ replaced by $\overline x,\overline y$, it holds that $(\pathspVar{1}(V),\overline\rho)$ is a genuine metric space (which is sensitive to the parameterization of paths).

\begin{corollary}\label{Cor:unique and parameterization}
The statement of Theorem~\ref{thm:universal map BV paths} holds with 
$\pathsreduced{1}$ replaced by $(\pathspVar{1},\overline\rho)$
and $\Phi$ replaced by the map $
\overline\Phi : \pathspVar{1}(V) \to \FSeriesOne{\overline V}$,
$x \mapsto \Lambda \circ \Sig(\overline x)$.
\end{corollary}
\begin{proof}
The above remarks imply that the map $x\mapsto \overline x$
injects $(\pathspVar{1}(V),\overline\rho)$ continuously into $\pathsreduced{1}(\overline V)$,
and the conclusion follows from Theorem~\ref{thm:universal map BV paths}. 
\end{proof}
A particularly important corollary is that that the robust signature moments characterize laws of stochastic processes.

\begin{corollary}\label{cor:BV_measure_sep}
Let $\mu,\nu$ be probability measures on $\pathspVar{1}(V)$. 
Then
\[
\mu = \nu \quad \text{ iff } \quad \EE_{X\sim \mu}\big[\overline\Phi(X)^m\big] = \EE_{Y\sim\nu}\big[\overline\Phi(Y)^m\big] \quad \forall m \geq 1\;.
\]
\end{corollary}
Corollary~\ref{cor:BV_measure_sep} is in stark contrast to classical signature moments, i.e.~it does not hold when $\overline\Phi(X),\overline\Phi(Y)$ are replaced by $\Sig(\overline X),\Sig(\overline Y)$; see Appendix~\ref{ex:example nonchar}.

\subsection{Lifting State Space Features to Pathspace Features}\label{subsec:state to path features}
So far we have focused on paths that evolve in a linear space.
In this section we briefly discuss how one can generalize the robust signature maps from Sections~\ref{subsec:univPaths} and~\ref{subsec:timepar BV} into a feature map for paths that evolve in general topological spaces.


\begin{example}
Let $\calX$ be the space of finite graphs\footnote{Say, directed and weighted graphs without multiple edges; however, any definition of a graph for which an adjacency matrix is well-defined is allowed.}
with vertex set $\{1,\ldots, N\}$.
Consider the feature map $\varphi: \calX \to E \coloneqq \RR^{N\times N}$ that maps a graph to its adjacency matrix.
Then $x \to (\Phi \circ \varphi)(x)$ is a natural candidate for a feature map for paths evolving in the space of graphs (an evolving network).
Other choices for $\varphi$ are of course possible and often more sensible, e.g.~those arising from kernels and we return to this in Section~\ref{subsec:kernel_learning}.
\end{example}

Even for paths that evolve in a linear space $\calX=\RR^d$, precomposition of $\Phi$ with a state-space non-linearity $\varphi$ can be beneficial in terms of efficiency since the state space non-linearity can make it easier to linearize non-linear functions.

\begin{definition}\label{def:normSigLift}
Let $\calX$ be a topological space and $\varphi : \calX\to V$.
Denote with $\pathspVar1_{x_0}(\calX)$ the subset of paths $x$ in $\cup_{T>0} C([0,T],\calX)$ such that $x(0)=x_0$ and $\varphi(x)=(\varphi(x(t)))_{t \in [0,T_x]} \in \pathspVar1(V)$.
Denote with $\pathsreduced{1}_{x_0}(\calX) = \pathspVar{1}_{x_0}(\calX) / \sim_t$ the quotient space\footnote{Equipped with the initial topology induced by the signature map $\Sig\circ\varphi : \pathsreduced{1}_{x_0}(\calX) \to \FSeriesOne{V}$.} under tree-like equivalence.
 \end{definition}
\begin{corollary}\label{cor:liftOfFeatures}
Suppose $\varphi:\calX \rightarrow V$ is injective and $\Lambda$ is a tensor normalization. 
The statement of Theorem~\ref{thm:universal map BV paths} holds with $\pathsreduced{1}(V)$ replaced by $\pathsreduced{1}_{x_0}(V)$
and $\Phi$ replaced by the map
\[
  \Phi^\varphi: \pathspVar{1}(\calX) \rightarrow \FSeriesOne{V}, \quad \Phi^\varphi \coloneqq \Phi \circ \varphi\;.
\]
Furthermore, the same statement holds when $\pathsreduced{1}_{x_0}(\calX)$ is replaced by $(\pathspVar{1}_{x_0}(\calX),\overline\rho)$ and $\Phi^\varphi$ is replaced by $\overline \Phi^\varphi\coloneqq\Lambda \circ \Sig (\overline{\varphi(x)})$.
\end{corollary}
\begin{proof}
Injectivity of $\varphi$ implies that for all $x,y \in \pathspVar{1}_{x_0}(\calX)$, $\varphi(x)\sim_t\varphi(y)$ if and only if $x \sim_t y$.
All the desired claims now follow from Theorem~\ref{thm:universal map BV paths} and Corollary~\ref{Cor:unique and parameterization}.
\end{proof}
Applied with $\calX = V$ and $\varphi=\id : \calX\to V$, this just recovers the feature map $\Phi^{\id}=\Phi$ for $\pathsreduced 1$.
Lifts of $\calX$-valued sequences were introduced for sequence kernels in~\citep{OK16}.

\section{Robust Features for Rough Paths}\label{sec:rough paths}
In this section we enlarge the domain of the feature map $\Phi$ (resp.~$\overline \Phi$) to paths of unbounded 1-variation.
This allows to include important classes of stochastic processes such as stochastic differential equations, semimartingales, Markov processes, and Gaussian processes.
The main obstacle is that the iterated integrals $\int dx^{\otimes m}$ can no longer be defined by Riemann--Stieltjes integration since the trajectories of $x$ are too irregular.
Rough paths theory provides a complete integration theory for a large class of processes; for the special case of (semi)martingales, the theory agrees with It{\^o} integration.
We focus here on the main ideas but give a detailed account and examples in Appendix~\ref{appendix:rough}.

Rough path theory provides a family of spaces $(\pathspVar{p})_{p \geq 1}$ called geometric $p$-rough paths over $V$.
These spaces satisfy the inclusions $\pathspVar{p} \subset \pathspVar{q}$ for $p \leq q$, and larger values of $p$ allow for paths with ``rougher'' trajectories.
For every $p\geq 1$, there exists a map
$
\pathspVar{p} \rightarrow \FSeriesOne{V}
$
which has the same properties as the signature map $\Sig$ from Section~\ref{sec: ordered moments}.
In fact, for $p=1$, $\pathspVar{1}$ is just the space of absolutely continuous paths and we recover the setting of Section~\ref{sec: ordered moments}.
Similar to before, we define the space $\pathsreduced{p} = \pathspVar{p}/\sim_t$, where $\sim_t$ denotes tree-like equivalence, which we equip with the topology induced by embedding $ \pathsreduced{p} \to \FSeriesOne{V}$.
The following is a significant generalisation of Theorem~\ref{thm:universal map BV paths} to large classes of processes, see Example~\ref{Ex:rough path lifts}. 
The proof follows however, analogous to the proof of Theorem~\ref{thm:universal map BV paths}. 
\begin{theorem}\label{thm:featureMapRPs}
Let $p \geq 1$ and $\Lambda : \FSeriesOne{V} \to \FSeriesOne{V}$ a tensor normalization.
The map
\begin{equation*}
  \Phi:\pathsreduced{p}\rightarrow \FSeriesOne{V}\;, \quad \Phi = \Lambda\circ\Sig 
\end{equation*}
  \begin{enumerate}[label=(\roman*)]
  \item is a continuous injection from $\pathsreduced{p}$ into a bounded subset of $\FSeriesOne{V}$,
  \item is universal to $\calF \coloneqq C_b(\pathsreduced{p},\RR)$ equipped with the strict
    topology,
  \item is characteristic to the space of finite regular Borel measures on $\pathsreduced{p}$.
  \end{enumerate}
  Moreover, $\forall \ell \in \FSeries{V}'$ the map $\mu \to \langle \ell, \EE_{X \sim \mu}[\Phi(X)]\rangle$ is $B$-robust.
  The same statement holds when $\pathsreduced{p}$ is replaced by $\pathspVar{p}$ and $\Phi$ is replaced by $\overline \Phi(\bfx)$, see Appendix~\ref{subsec:time par RPs}.
    \end{theorem}

\section{A Computable Metric for Laws of Stochastic Processes}\label{sec:kernel}
Sections~\ref{sec: ordered moments} and~\ref{sec:rough paths} provide us with a robust feature map $\Phi$ (resp.~$\overline \Phi$) such that its expected value $\EE_{X \sim \mu}[\Phi(X)]$ characterizes the law $\mu$ of the stochastic process $X$ [with the option to ignore parameterization]. 
In this section we use this to derive a metric on the space of laws of stochastic processes and study the topology it induces.

A natural distance for probability measures on any topological space $\calX$ is given by fixing a class of functions $\calG \subset \RR^\calX$ and define the \emph{maximum mean discrepancy} (MMD) between $\mu$ and $\nu$ as
\[
d(\mu,\nu) \coloneqq \sup_{f \in \calG} | \EE_{X\sim \mu}[f(X)] - \EE_{Y \sim \nu}[f(Y)] |\;.
\]
If $\calG$ is sufficiently large, then $d$ defines a metric between probability measures, see~\citep{rachev1991probability,muller1997integral}.
However, the supremum makes it hard to compute or estimate $d(\mu,\nu)$.
An insight from kernel learning is that if $\calG$ is the unit ball of a RKHS $(\calH,\kernel)$, then the reproducing property implies
\begin{equation}\label{eq:MMD_expectations}
d^2(\mu,\nu)=\EE[\kernel(X,X^\prime)]-2\EE[\kernel(X,Y)]+\EE[\kernel(Y,Y^\prime)]\;,
\end{equation}
where $X,X'\sim \mu$ and $Y,Y'\sim \nu$ are pairwise independent.
The above can be simply estimated from finite samples of $\mu$ and $\nu$ provided $\kernel$ is cheap to evaluate.

We are interested in the case when $\calX=\pathspVar{p}$ is the space of $p$-rough paths and $\mu$ and $\nu$ are the laws of stochastic processes $X=(X_t)$ and $Y=(Y_t)$.
The robust signature feature map $\Phi$ gives us a natural candidate to define a kernel $\kernel(x,y)\coloneqq\langle \Phi(x),\Phi(y) \rangle$.
A well-known disadvantage of MMDs is that they metrize a topology that is usually not comparable (i.e.~neither weaker nor stronger) to weak convergence, see~\citep{rachev1991probability, muller1997integral}. 
Progress on this question was made recently in~\citep{SGS18,simongabriel2020metrizing} in the case of locally compact spaces $\calX$ building on \citep{sriperumbudur2016optimal, sriperumbudur2011universality}. 
However, already the space of absolutely continuous paths $\pathspVar{1}$ is not locally compact.

\subsection{Kernel Learning}
\label{subsec:kernel_learning}
We briefly recall notation and terminology from kernel learning~\citep{cucker2002mathematical,scholkopf2002learning,berlinet2011reproducing}.
Throughout this section, let $\calX$ be a topological space.
Let us fix a feature map $\calX\to\RR^\calX$, $x\mapsto \kernel_x$,
and denote by $\calH_0\subset \RR^\calX$ the linear span of $\{\kernel_x:x\in \calX\}$.
Under the condition that $\kernel(x,y) \coloneqq \kernel_x(y)$ is a symmetric positive definite kernel (which we assume henceforth), the completion of $\calH_0$ under the inner product $\scal{\kernel_x,\kernel_y} \coloneqq \kernel(x,y)$ is the RKHS of $\kernel$.
We further fix a locally convex TVS $\calF\subset \RR^\calX$ containing $\calH_0$ for which the inclusion map $\calH_0 \hookrightarrow \calF$ is continuous.
      \begin{definition}[\citet{SGS18}]
        We say that the kernel $\kernel$ is
        \begin{itemize}
        \item universal to $\calF$ if the {kernel embedding} $\iota = \id :\calH_0 \hookrightarrow \calF, \;f\mapsto f$,
          is dense,
         \item characteristic to $\calF'$ if the {kernel mean embedding} $\mu:\calF^\prime
          \rightarrow \calH_0', \;D\mapsto D|_{\calH_0}$ is injective.
        \end{itemize}
        \end{definition}
        \begin{proposition}[\citet{SGS18}]\label{prop:kernel mean iff char}
The kernel embedding $\iota$ and the kernel mean embedding $\mu$ are transpose: $
\iota^\star = \mu$ and $\mu^\star= \iota$.
Moreover, $\kernel$ is universal to $\calF$ iff $\kernel$ is characteristic to $\calF^\prime$.
        \end{proposition}     
A common way to construct the feature map $x\mapsto\kernel_x$ is through another map $\Phi : \calX \to E$ taking values in an inner product space $(E,\scal{\cdot,\cdot})$.
Then one defines $\kernel_x(y) \coloneqq \kernel(x,y) \coloneqq \langle \Phi(x),\Phi(y) \rangle$.
We recall the important details of this construction in Appendix~\ref{appendix:kernels}.
As expected, universality and characteristicness of $\Phi$ and $\kernel$ are equivalent.
\begin{proposition}\label{prop:kPhiUnivChar}
Consider a map $\Phi:\calX \rightarrow E$ into an inner product space $(E,\langle \cdot,\cdot \rangle)$ such that $\kernel(x,y) = \scal{\Phi(x),\Phi(y)}$.
Then
        \begin{itemize}
        \item $\Phi$ is universal to $\calF$ iff $\kernel$ is universal to $\calF$,
        \item $\Phi$ is characteristic to $\calF^\prime$ iff $\kernel$ is characteristic to $\calF^\prime$.
       \end{itemize}
      \end{proposition}
\begin{proof}
This follows from Proposition~\ref{prop:kPhiUnivCharFull} in Appendix~\ref{appendix:kernels}.
\end{proof}

\subsection{The Robust Signature Kernel and its MMD}\label{subsec:kernel mean embedding}
The feature map $\Phi$ defined in Section~\ref{sec:rough paths} is universal and characteristic for (unparameterized) paths that evolve in a Banach space $V$.
It takes values in the feature space $\FSeries{V}$ that is itself a Banach space.
In particular, $\Phi$ can be applied to (unparameterized) paths that evolve in a Hilbert space $(H,\langle \cdot, \cdot \rangle)$.
In this case the feature space $\FSeries{H}$ is also a Hilbert space, see Section~\ref{sec:normFunc}.
This allows to take an inner product in the feature space which immediately gives rise to a kernel on $\pathsreduced{p}(H)$ that is also universal and characteristic.
\begin{theorem}\label{thm: signature kernel}
Let $p \geq 1$ and $\Normal:\FSeriesOne{H}\to\FSeriesOne{H}$ a tensor normalization.
Define
\[
  \kernel:\pathsreduced{p}(H)\times \pathsreduced{p}(H)\rightarrow \RR,\quad(x,y)\mapsto \langle \Phi(x),\Phi(y) \rangle
\]
where $\Phi=\Normal\circ \Sig$ is the normalized signature map.
We call $\kernel$ the robust signature kernel and denote with 
\[
  d_{\kernel}(\mu,\nu) \coloneqq\sup_{f \in \FSeries{H}:\Norm{f} \le 1} \norm{\int f(x) \mu(dx) -\int f(y) \nu(dy)}
\]
the associated MMD.
Then $\kernel$ is a bounded, continuous, positive definite kernel and 
\begin{enumerate}[label=(\roman*)]
\item $\kernel$ is universal to $C_b(\pathsreduced{p}(H),\RR)$ equipped with the strict topology,
\item $\kernel$ is characteristic to finite, signed Borel measures on $\pathsreduced{p}(H)$,
\item $d_{\kernel}$ is a metric on the space of finite, signed Borel measures on $\pathsreduced{p}(H)$, and
  \[
d^2_{\kernel}(\mu,\nu)=\EE[\kernel(X,X^\prime)]-2\EE[\kernel(X,Y)]+\EE[\kernel(Y,Y^\prime)]
  \]
  where $X,X'\sim \mu$ and $Y,Y'\sim \nu$ are pairwise independent,
  \item \label{point:topology comparison} for probability measures on $\pathsreduced{p}(H)$, convergence in $d_{\kernel}$ does not imply weak convergence.
Suppose now that $H$ is finite dimensional.
Then weak convergence implies convergence in $d_{\kernel}$.
Furthermore, if $\calM$ is a set of probability measures on $\pathsreduced{p}(H)$ which is compact under weak convergence, then $d_\kernel$ induces the topology of weak convergence on $\calM$.
\end{enumerate}
The same statement holds when $\pathsreduced{p}(H)$ is replaced by $(\pathspVar{p}(H),\overline\rho_p)$ and $\Phi$ in $\kernel(x,y)=\langle \Phi(x),\Phi(x) \rangle$ is replaced by $\overline \Phi(\bfx) \coloneqq\Lambda \circ \Sig (\overline \bfx)$.
\end{theorem}
\begin{proof}
All claims, except~\ref{point:topology comparison}, follow from combining Theorem~\ref{thm:featureMapRPs} and Proposition~\ref{prop:kPhiUnivChar}.
The proof of~\ref{point:topology comparison} is more involved and follows from Proposition~\ref{prop:weaker than weak} in Appendix~\ref{appendix:MMD for RPs}.
\end{proof}

\subsection{Lifting State Space Kernels to Pathspace Kernels}\label{sec: lifting kernels}
In Section~\ref{subsec:state to path features} we discussed how precomposing the feature map $\Phi$ with a state space non-linearity $\varphi$ can be beneficial. 
The same applies to the robust signature kernel $\kernel$.
\begin{proposition}\label{prop:lift_kernel_univ}
  Let $(H,\kappa)$ be a RKHS on $\calX$ such that $\varphi : \calX \to H$, $\varphi(x) \coloneqq\kappa(x,\cdot)$, is injective\footnote{$\varphi$ is injective whenever $\kappa$ is characteristic to, say, the space of finite, signed Borel measures on $\calX$.} and continuous.
  For $x \in C([0,T],\calX)$ denote by $x \mapsto \kappa_x \coloneqq(\varphi(x(t)))_t \in C([0,T],H)$ the lift of $x$ to a path evolving in $H$.
We call\footnote{Note that $\kernel$ in $\kernel(\kappa_x,\kappa_y)$ is the signature kernel on $\pathsreduced{1}(H)$ from Theorem~\ref{thm: signature kernel}.}
\[
  \kernel^{\kappa}:\pathsreduced{1}(\calX)\times \pathsreduced{1}(\calX) \rightarrow \RR, \quad (x,y)\mapsto\kernel(\kappa_x,\kappa_y)
\]
the robust signature lift of the kernel $\kappa: \calX \times \calX \rightarrow \RR $.
Then $\kernel^{\kappa}$ is a bounded, continuous, positive definite kernel and, using the notation from Definition~\ref{def:normSigLift}, for any $x_0\in\calX$
\begin{enumerate}[label=(\roman*)]
\item $\kernel^{\kappa}$ is universal to $C_b(\pathsreduced{1}_{x_0}(\calX),\RR)$ equipped with the strict topology,
\item $\kernel^{\kappa}$ is characteristic to finite, signed Borel measures on $\pathsreduced{1}_{x_0}(\calX)$,
\item $d_{\kernel^\kappa}$ is a metric on the space of finite, signed Borel measures on $\pathsreduced{1}_{x_0}(\calX)$.
\end{enumerate}
The same statement holds when $\pathsreduced{1}_{x_0}(\calX)$ is replaced by $(\pathspVar{1}_{x_0}(\calX),\overline\rho)$ and $\Phi$ (implicit in the definition of $\kernel$) is replaced by $\overline\Phi(x) \coloneqq\Lambda \circ \Sig (\overline x)$.
\end{proposition}
If $\calX$ is itself a Hilbert space, the choice $\kappa(\cdot, \cdot ) = \langle \cdot, \cdot \rangle$ corresponds to $\varphi = \id$ and recovers the signature kernel $\kernel$ from Theorem \ref{thm: signature kernel}, $\kernel^{\langle \cdot, \cdot \rangle} = \kernel$.
\subsection{Computing the Signature Kernel}\label{subsec:comp_sig_kernel}
We now turn to the computational aspects of the robust signature kernel $\kernel^\kappa$ and associated MMD.
In practice, we are not given a path sample $x=(x(t))_{t \in [0,T]}$, but a sequence sample $x^\pi=(x(t_i))_{i}$ that measures the values of $x$ at times $t_i$ in a partition $\pi=\{0 \le t_1 < \cdots \le T\}$ of $[0,T]$.
This partition $\pi$ might even change from path to path in a given dataset.
Hence, instead of working with probability measure on paths we work with probability measures on sequences.
Below we define the discrete robust signature (resp.~kernel resp.~MMD) on the domain of sequences and then quantify how far these are from the robust signature $\Phi^\varphi$ (resp.~kernel $\kernel^\kappa$ resp.~MMD $d_{\kernel^\kappa}$) that is defined on the domain of paths.
Further, we show that the algorithms of~\citet{OK16} can be used to compute these quantities efficiently.

For the rest of this section, we fix an integer $M\geq 1$, a RKHS $(H,\kappa)$ on $\calX$, and a $K$-Lipschitz function $\psi$ as in Corollary~\ref{cor:tensor_normal} with associated tensor normalization $\Lambda$.

\begin{definition}\label{def:BVApproximations}
  Denote by $\calX^+ = \bigcup_{\ell\ge 0} \calX^\ell$ the space of sequences in $\calX$ of arbitrary length.
  Define
  \begin{equation}
    \kernelTrunc{M}^{\kappa,+}:  \calX^+\times  \calX^+\rightarrow \RR\;, \quad (x,y) \mapsto \scal{\Phi_{M}^{\linear}(\kappa_{x}),\Phi_{M}^{\linear}(\kappa_{y})}
  \end{equation} 
  where, for $x=(x_0,\ldots, x_\ell)$, $\kappa_x$ is the $H$-valued sequence $(\kappa(x_0,\cdot),\ldots, \kappa(x_\ell,\cdot))$ and\footnote{We use the shorthand notation $1+\kappa(x_i,\cdot)-\kappa(x_{i-1},\cdot)\coloneqq(1,\kappa(x_i,\cdot)-\kappa(x_{i-1},\cdot), 0,\ldots,0) \in  \TSeriesTruncated{H}{M}$ and $\prod_{i=1}^\ell x_i$ denotes the product $x_1\otimes x_2\otimes\cdots \otimes x_\ell$ in the tensor algebra, see Appendix~\ref{sec:tensor algebra}.}
 
  \begin{equation}
    \Phi_{M}^{\linear}: \calX^+ \to \TSeriesTruncated{H}{M}, \quad   \Phi_{M}^{\linear}= \Lambda \Big( \prod_{i=1}^\ell  (1+\kappa(x_i,\cdot)-\kappa(x_{i-1},\cdot))\Big)\;.
  \end{equation}
\end{definition}

As we show in Proposition~\ref{prop:kernel_discrete} of Appendix~\ref{appendix:sig_kernels},
$\kernelTrunc{M}^{\kappa,+}:  \calX^+\times  \calX^+\rightarrow \RR$ is a bounded, positive semidefinite kernel for $\calX^+$.
The reason why $\kernelTrunc{M}^{\kappa,+}$ is of practical interest is that it approximates $\kernel^{\kappa}$
with an explicit convergence rate.
More precisely.
for $x,y \in \pathspVar{1}(\calX)$ and partitions $\pi \subset [0,T_x],\pi'\subset[0,T_y]$,
\begin{equation*}
\big| \kernel^{\kappa}(x,y)- \kernelTrunc{M}^{\kappa,+}(x^\pi, y^{\pi'}) \big|
\leq
C
\big(
\max_{t_i\in \pi}\|\kappa_x\|^{1/2}_{1\Var;[t_{i-1},t_{i}]} 
+ \max_{t_i\in \pi'}\|\kappa_y\|^{1/2}_{1\Var;[t_{i-1},t_{i}]}
\big)\;,
\end{equation*}
where $x^\pi$ (resp. $y^{\pi'}$) denote the sequences given by sampling $x$, $y$ along $\pi$ (resp.~$\pi'$),
and where $C$ depends only on $\|\psi\|_\infty, K,\|\kappa_x\|_{1\Var}$, and $\|\kappa_y\|_{1\Var}$.
A similar approximation result holds for the associated MMDs $d_{\kernel^\kappa}$ and $d_{\kernelTrunc{M}^{\kappa,+}}$.
See again Proposition~\ref{prop:kernel_discrete} for a precise statement.

The following result provides a bound on the computation cost of $\kernelTrunc{M}^{\kappa,+}$.

\begin{proposition}\label{prop:kernel_complexity}
The kernel $\kernelTrunc{M}^{\kappa,+}(x^\pi, y^{\pi'})$ can be evaluated in $O((|\pi|^2+ |\pi'|^2)(c + M) + q)$ time and $O(|\pi|^2+ |\pi'|^2 + M + r)$ memory, where $|\pi|$ is the number of time points in $\pi$, $c$ is the cost of one evaluation of the kernel $\kappa$, and $q$ and $r$ are the total time and memory costs respectively of a single evaluation of $\psi$ and of finding the unique non-negative root of a polynomial $P(\lambda) = \sum_{m=0}^{M} a_{m}\lambda^{2m}$ with $a_0 \leq 0 \leq a_1,\ldots, a_M$.
\end{proposition}

\begin{proof}
Denote $\Sig_{M}^{\linear}(h) \coloneqq \prod_{i=1}^\ell (1+h_i - h_{i-1})$.
The sequence of norms
$\|(\Sig_{M}^{\linear}(\kappa_{x^\pi})^m\|_{H^{\otimes m}}$, $m=0,\ldots, M$, can be computed in $O(|\pi|^2(c + M))$ time and $O(|\pi|^2 + M)$ memory due to \citep[Alg.~3]{OK16}.
In the same way, the sequence of inner products
$(\scal{\Sig_{M}^{\linear}(\kappa_{x^\pi})^m,\Sig_{M}^{\linear}(\kappa_{y^{\pi'}})^m}_{H^{\otimes m}})_{m=0}^M$ can be computed in $O(|\pi||\pi'|(c + M))$ time and $O(|\pi||\pi'| + M)$ memory.
We now claim that, for $\bfs,\bft\in \TSeriesTruncated{H}{M}$, once the three
sequences $(\|\bfs^m\|)_{m =0}^M$, $(\|\bft^m\|)_{m =0}^M$,
and $(\scal{\bfs^m,\bft^m})_{m =0}^M$
are given,
we can compute the normalized inner product $\scal{\Lambda(\bfs),\Lambda(\bft)}$
in $O(q)$ time and $O(r)$ memory,
which will finish the proof.
Indeed, let $\lambda(\bft)\geq 0$ be as in Proposition~\ref{prop:LipschitzBoundPhi}.
Then $\scal{\Lambda(\bfs),\Lambda(\bft)} = \sum_{m=0}^M \lambda(\bfs)^m\lambda(\bft)^m\scal{\bfs^m,\bft^m}$,
so it suffices to compute $\lambda(\bfs)$ and $\lambda(\bft)$.
These can be computed in the claimed time and memory since $\lambda(\bft)$ is the unique non-negative root of the polynomial
$
P(\lambda)  \coloneqq \|\delta_{\lambda} \bft \|^2 -  \psi(\|\bft\|) = \Big(\sum_{m=0}^M\lambda^{2m} \|\bft^{m}\|^2 \Big) -  \psi(\|\bft\|)$.
\end{proof}
The quadratic complexity in the sequence length $|\pi|$ can be reduced to linear complexity by using the simultaneous dynamic programming and low-rank algorithm from \citep[Algorithm 5, Section 6]{OK16}; however, although the low-rank assumption works well in practice, it is unclear how to justify it theoretically in the context of paths and signatures. 
\section{Application: Two-sample Tests for Stochastic Processes}\label{sec:twosample}
We apply the MMD metric developed in the previous sections to the two-sample testing problem
  \begin{align}
    H_0: \mu = \nu \text{ against the alternative } H_1: \mu \neq \nu
  \end{align}
where $\mu\coloneqq \operatorname{Law}(X) \coloneqq\PP \circ X^{-1}$ and $\nu \coloneqq \operatorname{Law}(Y) \coloneqq\PP \circ Y^{-1}$ denote the laws of the stochastic processes $X$ and $Y$.
We compare the resulting MMD test against classical multi-variate tests as well as MMD tests using classical kernels for vector-valued data. 
\subsection{Test statistics} 
Denote with $X=(X(t_1), \ldots,X(t_L)), Y=(Y(s_1), \ldots, Y(s_L))$ two independent, discrete-time stochastic processes that evolve in $\RR^d$. 
We are given $m$ (resp.~$n$) i.i.d.~samples from $X$ (resp.~$Y$), which we denote with $\bX=\{X_1,\ldots,X_m\}$ and $\bY=\{Y_1,\ldots,Y_n\}$.
That is, each $X_i$ is a sequence of $L$ vectors in $\RR^d$. 
The probability of falsely rejecting the null is called the \emph{Type I error}, and similarly the probability of falsely accepting
the null is called the \emph{Type II error}.
If the Type I error can be bounded from above, uniformly over all $\PP$ under which
$X$ and $Y$ are independent, by a constant $\alpha$, then we say that the test is of \emph{significance level} $\alpha$.
Given a statistic $T(\bX,\bY)$ one can construct a simple permutation test: under the null $H_0$ one can sample from the distribution of $T(\bX,\bY)$ by sampling uniformly at random from the $(m+n)!$ permutations $\pi$ of $\{1,\ldots,m+n\}$ and evaluate
\begin{align*}
 T(Z_{\pi(1)},\ldots,Z_{\pi(m)},Z_{\pi(m+1)},\ldots,Z_{\pi(m+n)})
\end{align*}
where $\bZ=(X_1,\ldots,X_m,Y_1,\ldots,Y_n)$.
In particular, this allows to estimate the $95\%$ quantile of $T(\bX,\bY)$ under the null $H_0$. 
Hence, if $T(\bX,\bY)$ is not in this $95\%$ quantile, we can reject the null with significance level $5\%$.
For background on testing and (permutation) test statistics we refer to
\citep[Chapter 15]{lehmann2006testing}. 
\paragraph{MMD statistics.}
Using the sequence signature kernel $\kernel^{\kappa,+}_M$ from Definition~\ref{def:BVApproximations}, we can use the MMD
  \begin{align}
    \label{d_k}
    d_{\kernel^{\kappa,+}_M}(\mu,\nu)^2 =\EE[\kernel^{\kappa,+}_M(X,X')]-2\EE[\kernel^{\kappa,+}_M(X,Y)]+\EE[\kernel^{\kappa,+}_M(Y,Y')]
  \end{align}
where $X,X'$ (resp. $Y,Y'$) denote pairwise independent sequences sampled from $\mu$ (resp.~$\nu$).
The unbiased (quadratic complexity) estimator from \citep{gretton2009fast,gretton2012kernel} turns \eqref{d_k} into a test statistic $T(\bX,\bY)$ and combining this with the results of Section~\ref{sec:kernel} allows to efficiently estimate \eqref{d_k}. 
Moreover, we can do the same by replacing the robust sequence signature kernel $\kernel^{\kappa,+}_M$ in \eqref{d_k} with any kernel  $\kappa$ for vector-valued data by flattening sequence samples; e.g.~if $X$ evolves in the state space $\RR^d$, we identify the sequence of vectors $(X(t_1),\ldots,X(t_L))$ as a single vector in $\RR^{d\cdot L}$. 
  The choices for the kernel $\kappa$ we used (both in the signature MMD and for the vector kernel MMD of the flattened sequence) were the linear kernel $\kappa(x,y)=\scal{x,y}$,
  the RBF kernel $\kappa(x,y) = \exp(-\|x-y\|^2/(2\sigma^2))$,
  and the Laplace kernel $\kappa(x,y) = \exp(-\|x-y\|/\sigma)$ for a parameter $\sigma>0$.
  Kernels depend on parameters:
  for the signature MMD these parameters are the truncation level $M$, the parameters $a,C$ that determine the tensor normalization (Example \ref{ex:psi}), and the parameters of the kernel $\kappa$;
  for the classical kernel choices of $\kappa$ the only parameter is the length scale $\sigma$.

\paragraph{Classic multi-variate statistics.}
By flattening a sequence observations into a long vector (as described above), one can also use classical multi-variate test statistics $T(\bX,\bY)$.
Many of these classical tests like Biau-Gy\"orfi scale poorly and the high-dimensionality of our testing problem makes it impossible to use them; see \citep{gretton2012kernel} for more on this.
As in \citep{gretton2012kernel} we settled for the Friedman--Rafsky multivariate Wald--Wolfowitz test statistic from \citep{friedman1979multivariate}, (Wolf), and the Hotelling t-test statistic, (Hotelling), since they are still computable in moderately high-dimension.
However, even here, computing these tests requires building a graph for (Wolf) and matrix inversion for (Hotelling).
This makes them scale at least cubically in $dL$ and thus are only applicable if both, the dimension of the state space $d$ and the sequence length $L$, are small.   
\subsection{Data}\label{sec:data}
As data sources we use the multi-variate time series datasets from UAE \& UCR, see \citep{bagnall2018uea} for a detailed description. 
This archive consists of 30 time series datasets from various domains, including speech, motion, handwriting, robotics, etc.
The state space $\RR^d$ of these datasets varies between low-dimensional $d=2$ (e.g.~Pendigits) to high-dimensional $d={1345}$ (e.g.~DuckDuckGeese).

For each time-series dataset we define $\mu$ (resp.~$\nu$) as follows: under the null $H_0$ a sample from $\mu$ (resp.~$\nu$) is generated by choosing uniformly at random a time-series $W=(W(t_1),\ldots,W(t_{L_0}))$ from this dataset and double the sequence length by concatenation; that is $X=(W(t_1),\ldots,W(t_{L_0}),W(t_1),\ldots, W(t_{L_0}))$.
Under the alternative $H_1$, a sample from $\mu$ is generated as above, but a sample from $\nu$ is generated by concatenating in reverse, $Y=(W(t_1),\ldots,W(t_{L_0}),W(t_{L_0}),\ldots, W(t_1))$.
Despite its simplicity, we will see that this problem is challenging for many statistics, since a statistic needs to capture the order structure to achieve a high test power. 
In our experiments, we varied the length $L_0 \in \{ 10, 100, 200\}$ by truncating after the first $L_0$ entries of the time-series in the original dataset. 
Further, we restricted to the case of balanced samples, $m=n$, and varied the number of samples $m \in \{30,70,200\}$.
This allows to also compare the sensitivity of the different statistics to length and sample size.

\subsection{Implementation Details}\label{subsec:implementation}
As base for our signature MMD computations we use Csaba Toth's KSig library\footnote{ \url{https://github.com/tgcsaba/KSig}} that provides efficient GPU implementations of signature kernels in Python.
As base for general kernel MMD two-sample permutation test, we built on code from the DS3 summer school course ``Representing and comparing probabilities with kernels'' provided by H.~Strathmann and D.~Sutherland\footnote{\url{https://github.com/karlnapf/ds3_kernel_testing}} but we rewrote much of it to make it compatible with the KSig package.  
To find the non-negative root of the polynomial $P(\lambda)$ from Proposition~\ref{prop:kernel_complexity} we used an optimization of Brent's method~\citep[Ch.~3-4]{Brent73} implemented as \texttt{optimize.brentq} in the \texttt{SciPy} package~\citep{python}.

\subsection{Discussion of results}
For each test statistic we used the permutation test with significance level $5\%$ and repeated the following for each dataset $20$ times:
\begin{enumerate*}[label=(\roman*)]
\item 
  Generate samples $\bX$,$\bY$ under the null $H_0$ and record the number of false rejections of the null $H_0$ (\emph{Type I error}).
  \item 
    Generate samples $\bX$,$\bY$ under the alternative $H_1$ and record the number of false acceptances of the null $H_0$ (\emph{Type II error}). 
\end{enumerate*}
The power of the test---the probability of rejecting the null $H_0$ when $H_1$ is true---then determines which test performs best; the closer to $1$ the better. 
As mentioned above for classical statistics the above was computationally only possible when the effective dimension is small hence these are only included in the tables when feasible. 

Table~\ref{tab:Type1L10m30} and Table~\ref{tab:Type2L10m30} show the empirical Type I error and $1$ minus the Type II error (i.e.~the power of the test) for the different statistics for sequence length $L_0=10$ and $m=30$. 
Table entries are rounded to two decimal digits and the most powerful test for each dataset is highlighted in bold.
The remaining results for the Type II error with varying sequence length $L_0$ and sample size $m$ are shown in Tables~\ref{tab:Type2L100m30}, \ref{tab:Type2L200m30}, and \ref{tab:Type2L200m200}.
For brevity we do not include tables with the Type I error for these other regimes since, by construction of the permutation test, these errors converge to a number less than $5\%$ as $m$ increases (that is, they all look very similar to Table~\ref{tab:Type1L10m30}, just more entries are less or equal to $0.05$).

The experiments show that, as one expects, the linear (inner product) kernel MMD does poorly overall and generally fails to distinguish two different distributions.
The standard MMDs with RBF and Laplace kernel outperform the linear kernel, but in turn rarely reach the test power of any of the signature MMDs.
As the length increases and thus samples carry more distinguishable structure,
the performance of all MMDs improves.
However, qualitatively the same behaviour remains.
The same is true when the number of samples $m$ increases; see Tables~\ref{tab:Type2L10m30}-\ref{tab:Type2L200m200}.
Among the classical statistics, the Hotelling test statistic is not competitive, even for small sequence length.
The Friedman-Rafsky multivariate version of the Wald-Wolfowitz test statistic performs well overall and is only beaten by the linear Signature MMD.
However, this test comes at a huge computational cost even for short sequences since it involves computing a minimal spanning tree.
Similarly Hotelling involves a matrix inversion.
This is the reason why can report these two, Hotelling and Wolf, only for the shortest sequence length $L_0=10$ and even there only for time-series that evolve in relatively low-dimensional state spaces (these experiments took several days on a multicore machine,\footnote{All experiments were run on a machine with 36 cores (2 x 18 Core with HT), 2.6/3.9 GHz, 1536GB RAM.} as opposed for the MMD experiments which took about an hour).  

Overall, these results show that signature MMDs efficiently exploit the sequence structure whereas classical MMDs need much more data to reach a similar test power.
A more surprising observation is that for the signature MMDs, the linear kernel does better than non-linear kernels.
This is in stark contrast to supervised time-series classification where non-linear signature kernels strongly outperform the linear signature kernel
as in~\citep{2019arXiv190608215T}.
We believe the reason is that the additional parameter selection for the non-linear kernel introduces additional variance, but unlike for the non-signature MMDs, the linear signature MMD is already expressive enough so that the additional non-linearities only result in more variance.
This trade-off might change considerably if one switches to other data resp.~hypothesis tests.

Related to this is also the question of kernel parameter choice.
All our tables show the results when the approach in \citep{fukumizu2009kernel} is applied for the kernel parameter choice.
We also experimented with the median heuristic as well as \citep{sutherland2016generative} that uses a part of the data for training to directly minimize the Type II error.
The median heuristic was not competitive and the price for learning the parameters is that the training data can not be used for the test itself and the choice of the train-test split is itself a heuristic.
In our experiments this choice had big effects from dataset to dataset so that ultimately \citep{fukumizu2009kernel} yielded overall the best and consistent performance for all MMDs. 
Nevertheless, improving the parameter choice for MMD tests is an ongoing research effort and any progress can potentially improve the results of all MMDs. 
The parameter selection resulted in a value for $M$ that was less or equal to $5$ for all signature MMDs.
Furthermore, for the signature MMDs build on top of a non-linear kernel $\kappa$, the value of $M$ was typically even lower, often equal to $2$.
This is in line with the discussion in Section~\ref{subsec:state to path features} and Section~\ref{sec: lifting kernels} that state-space non-linearities help to linearize so that a few iterated integrals suffice. 
In contrast, for the signature MMD build on top of the linear kernel, large values of $M$ could be in principle required. 
However, the price of using iterated integrals $\int dX^{\otimes M}$ for large values of $M$ is that their variance increases since we only have access to an empirical measure that approximates their expectation.
This explains why for the linear signature MMD, the optimal parameter $M$ is typically higher than for signature MMDs with non-linear $\kappa$ but still relatively low. 
The other parameters of the signature MMD $a,C$ determine the tensor normalization, Example \ref{ex:psi}.
For simplicity we kept both fixed ($a=1, C=10^3$)
since the theoretical guarantees apply to any choice and the signature MMD already outperform the baselines.
However, learning them could potentially improve the signature MMD results further. 
The parameter for the length scale $\sigma$ varied over several orders of magnitude from dataset to dataset, for both signature MMDs and the other MMDs. 

Finally, we note that one of the biggest strengths of the signature MMD is that it applies to sequences of different length.
None of the other statistics that identify a sequence as a long vector can be used here since the vectors that result from flattening the sequence would not have the same dimension. 
To test how much this influences the signature MMDs, we drew for each sequence sample a random integer $k$ between $0$ and $0.3 \cdot L_0$ and then deleted uniformly at random $k$ sequence entries.
The result is that each datasets contains sequences of varying length.
The signature MMDs can be applied without modifications and we repeated the above experiments.
The performance slightly deteriorates but remains very close to the signature MMD results given in the tables, hence we omit it for brevity.  

To sum up, we recommend to try both, the linear signature MMD and RBF signature MMD, for testing no matter what regime (small/large samples, small/large length, small/large state space).

\begin{table}
  \centering
  \scalebox{\tablescale}{
  \begin{tabular}{lccccccccc}
    \toprule
    & \multicolumn{3}{c}{MMDs} & \multicolumn{3}{c}{Signature MMDs} & \multicolumn{2}{c}{Classic Statistics} 
    \\\cmidrule(lr){2-4}\cmidrule(lr){5-7}\cmidrule(lr){8-9}
    & Linear  & RBF & Laplace & Linear & RBF & Laplace & Hotelling & Wolf\\\midrule
ArticularyWordRecognition                       & 0.1    & 0.05 & 0.1     & 0.05                & 0                & 0   &0.05&0.2        \\
AtrialFibrillation                              & 0.1    & 0.05 & 0.05    & 0                   & 0                & 0.1  &0.1&0.1       \\
BasicMotions                                    & 0.05   & 0    & 0.1     & 0.05                & 0.1              & 0.1  &0.08&0     \\
CharacterTrajectories                           & 0      & 0.1  & 0.05    & 0                   & 0.1              & 0.05 &0.1&0.15 \\
Cricket                                         & 0.05   & 0    & 0.1     & 0.1                 & 0.05             & 0 &0.1&0.05      \\
DuckDuckGeese                                   & 0.15   & 0.15 & 0.05    & 0.05                & 0                & 0.1 &-&-       \\
ERing                                           & 0      & 0.05 & 0.1     & 0.05                & 0                & 0 &0.05&0.1         \\
EigenWorms                                      & 0.05   & 0.05 & 0       & 0                   & 0.05             & 0.1&0&0       \\
Epilepsy                                        & 0.1    & 0.05 & 0.1     & 0.1                 & 0.05             & 0  &0&0.15            \\
EthanolConcentration                            & 0      & 0    & 0.05    & 0.1                 & 0                & 0   &0.05&0          \\
FaceDetection                                   & 0.05   & 0    & 0.05    & 0.05                & 0                & 0.1 &-&-           \\
FingerMovements                                 & 0      & 0.1  & 0       & 0.1                 & 0                & 0.05 &-&-          \\
HandMovementDirection                           & 0      & 0.05 & 0       & 0.05                & 0                & 0.05 &0&0.05         \\
Handwriting                                     & 0.15   & 0    & 0.05    & 0                   & 0.1              & 0.1  &0.1&0.05         \\
Heartbeat                                       & 0      & 0.05 & 0.15    & 0                   & 0.1              & 0.05  &-&-        \\
InsectWingbeat                                  & 0.05   & 0.25 & 0.05    & 0.1                 & 0                & 0     &-&-        \\
JapaneseVowels                                  & 0.25   & 0    & 0.05    & 0                   & 0                & 0.05 &0.05 &0.1      \\
LSST                                            & 0.05   & 0.1  & 0.15    & 0.05                & 0.1              & 0.05  &0&0.15        \\
Libras                                          & 0.1    & 0.05 & 0.1     & 0                   & 0.1              & 0.05  &0.1&0.1          \\
MotorImagery                                    & 0.05   & 0.1  & 0.1     & 0.1                 & 0.05             & 0.05  &-&-               \\
NATOPS                                          & 0.05   & 0.05 & 0.1     & 0                   & 0.05             & 0     &-&-               \\
PEMS-SF                                         & 0.05   & 0.05 & 0       & 0.05                & 0.05             & 0.05  &-&-               \\
PenDigits                                       & 0.05   & 0    & 0.15    & 0                   & 0                & 0.1   &0.05&0.1               \\
PhonemeSpectra                                  & 0.05   & 0.05 & 0.15    & 0.1                 & 0.05             & 0.05  &0.05&0        \\
RacketSports                                    & 0      & 0.05 & 0.05    & 0.2                 & 0                & 0.1   &0.1&0         \\
SelfRegulationSCP1                              & 0.15   & 0.05 & 0.1     & 0                   & 0.1              & 0     &0.05&0.05           \\
SelfRegulationSCP2                              & 0.05   & 0.05 & 0.05    & 0                   & 0.05             & 0.15  &0.05&0.1       \\
SpokenArabicDigits                              & 0      & 0.15 & 0.05    & 0.15                & 0.1              & 0    &0&0.05      \\
StandWalkJump                                   & 0.05   & 0.1  & 0       & 0.1                 & 0.15             & 0.1  &0&0.15         \\
UWaveGestureLibrary                             & 0.05   & 0    & 0.05    & 0.05                & 0                & 0.05 &0.05&0.05        \\
\bottomrule
\end{tabular}
}
\caption{Type I error for sequence length $L_0=10$ and $m=n=30$ samples.}
\label{tab:Type1L10m30}
\end{table}

\begin{table}[H]
  \centering
  \scalebox{\tablescale}{
  \begin{tabular}{lccccccccc}
    \toprule
    & \multicolumn{3}{c}{MMDs} & \multicolumn{3}{c}{Signature MMDs} & \multicolumn{2}{c}{Classic Statistics} 
    \\\cmidrule(lr){2-4}\cmidrule(lr){5-7}\cmidrule(lr){8-9}
    & Linear  & RBF & Laplace & Linear & RBF & Laplace & Hotelling & Wolf \\\midrule
ArticularyWordRecognition &	0 &	0 &	0.2 &	\ubold{0.55}	& 0.25& 	0.3& 	0.05	& 0.3 \\
AtrialFibrillation  & 0.2                                       & 0.2   & 0.15         & 0.7                 & 0.65             & \ubold{0.85}                  & 0.15                   & 0.8             \\
 BasicMotions                                                                                  & 0.2                                       & 0.45                                                                                                            & 0.45         & \ubold{1}                   & \ubold{1}                & \ubold{1}                    & 0.4                  &      \ubold{1}                    \\
CharacterTrajectories                                                                         & 0.3                                       & \ubold{1}   & \ubold{1}            & \ubold{1}                   & \ubold{1}               & \ubold{1}                    & 0.2                  & \ubold{1}                          \\
Cricket                                                                                       & 0.4                                       & 0.9                                                                                                             & 0.95         & \ubold{1}                   & \ubold{1}                & \ubold{1}              & \ubold{1}                     & \ubold{1}                          \\
DuckDuckGeese         & 0.05                                      & 0.05                                                                                                            & 0            & \ubold{0.9}                 & 0.45             & 0.2                  & - & - \\
ERing                                                                 & 0.05                              &0.45    &                        0.3          & \ubold{1}                   & 0.9              & 0.8                  & 0.3                  & 0.95                          \\
EigenWorms                                                                                    & 0                                         & 0.2                                                                                                             & {0.35}         & \ubold{1}                   & 0.15             & 0.15                 & 0.1                   & 0.15                         \\
Epilepsy                                                                                      & 0.1                                       & 0.15                                                                                                            & 0.1          & {0.35}                & 0.1              & 0.1                  & 0.1                  & \ubold{0.7}                          \\
EthanolConcentration                                                                          & 0.05                                      & 0.05                                                                                                            & 0.1          & \ubold{0.3}                & 0                & 0.05                 & 0                  & \ubold{0.3}                        \\
FaceDetection                                                                                 & 0                                         & \ubold{0.75}                                                                                                            & 0.15         & 0.2                 & 0.05             & 0.05                 & - &  - \\
FingerMovements                                           & 0.3                                       & 0.65                                                                                                            & 0.5          & \ubold{1}                   & 0.15             & 0.05                 & -                     & -                          \\
HandMovementDirection                                                                         & 0.05                                      & 0.05                                                                                                            & 0.05         & \ubold{0.95}                & 0.15             & 0.1                  & 0.05                     & 0.3                          \\
    Handwriting              & \ubold{1}                                       &  \ubold{1}                                                                                                               &  \ubold{1}           &  \ubold{1}                   &  \ubold{1}                &  \ubold{1}                    & 0.8                  & \ubold{1}                         \\
    Heartbeat                    & 0.05                                      & 0                                                                                                               & 0.1          & \ubold{0.25}                & 0                & 0.05                 & - & -                          \\
InsectWingbeat                                                                                & 0.05                                      & 0                                                                                                               & 0.05         & 0.1                 & \ubold{0.2}              & 0.1                  & - & - \\
    JapaneseVowels                                                                                &\ubold{1}                                         &  \ubold{1}                                                                                                               &  \ubold{1}   &  \ubold{1}                   &  \ubold{1}                &  \ubold{1}                   &  \ubold{1}                    & \ubold{1}                          \\
    LSST                                                                                          & 0.1                                       & 0.1                                                                                                             & 0.05         &  {0.35}               & 0.05             & 0.05                 & 0.14                 & \ubold{1}                          \\
    Libras                                                                                        & 0.05                                      & 0.05                                                                                                            & 0.2          & 0.15                &  {0.25}            & 0.15                 & 0.2                   & \ubold{0.95}                          \\
    MotorImagery                                                                                  & 0.05                                      & 0.05                                                                                                            & 0.05         &  \ubold{1}                  & 0.7              & 0.35                 &  - & - \\
    NATOPS                                                                                        & 0.1                                       & 0.15                                                                                                            & 0.15         &  \ubold{0.65}                & 0.15             & 0.1                  & -                  & -                          \\
    PEMS-SF                                                                                       & 0.05                                      & 0.1                                                                                                             & 0.05         &  \ubold{0.95}                & 0.25             & 0.25                 &  - & - \\
PenDigits                  & \ubold{1}              & \ubold{1}  & \ubold{1}            & \ubold{1}    & \ubold{1}          & \ubold{1}   & \ubold{1}               & \ubold{1}                          \\
    PhonemeSpectra                                                                                & 0.1                                       & 0.1                                                                                                             & 0            &  {0.3}                & 0.15             & 0.05                 & 0                     & \ubold{0.95}                          \\
    RacketSports                                                                                  & 0.3                                       & 0.25                                                                                                            & 0.3          &  \ubold{0.6}                 & 0.25             & 0.15                 & 0.24                 & 0.55                          \\
    SelfRegulationSCP1                                                                            & 0.15                                      & 0.1                                                                                                             & 0.1          & \ubold{1}                   & \ubold{1}               & \ubold{1}                    & 0.05                  & \ubold{1}                          \\
    SelfRegulationSCP2                                                                            & 0.05                                      & 0.1                                                                                                             & 0.05         & \ubold{1}                   & \ubold{1}               & \ubold{1}                   & 0.05                   & \ubold{1}                        \\
    SpokenArabicDigits                                                                            & 0.15                                      & 0.2                                                                                                             & 0.2          &  \ubold{1}                & 0.1              & 0.05                 & {0.35}                  & \ubold{1}                        \\
    StandWalkJump                                                                                 & 0.05                                      & 0.35                                                                                                            & 0.35         & 0.55                &  \ubold{0.9}             & 0.85                 & 0.35                   & 0.7                        \\
UWaveGestureLibrary                                                                           & 0                                         & 0                                                                                                               & 0.1          & \ubold{0.2}                 & 0.05            & 0.05                 & 0.05                  & 0                        \\
\bottomrule
\end{tabular}
}
\caption{1 minus the Type II error for sequence length $L_0=10$ and $m=n=30$ samples.}
\label{tab:Type2L10m30}
\end{table}

\begin{table}[H]
\centering
\scalebox{0.55}{
\begin{tabular}{lcccccccc}
\toprule
 & \multicolumn{3}{c}{MMDs} & \multicolumn{3}{c}{Signature MMDs} & \multicolumn{2}{c}{Classic Statistics} 
 \\\cmidrule(lr){2-4}\cmidrule(lr){5-7}\cmidrule(lr){8-9}
  & Linear  & RBF & Laplace & Linear & RBF & Laplace &  Hotelling & Wolf  \\\midrule
ArticularyWordRecognition                       & 0      & 0.1  & 0.35    & \ubold{1}             & 0.65             & 0.85             &0.15 &\ubold{1}    \\
AtrialFibrillation                              & 0.1    & 0.25 & 0.4     & \ubold{1}             & \ubold{1}     & \ubold{1}              & 0.6  & \ubold{1}\\
  BasicMotions                                    & 0.9    & \ubold{1} & 0.95    & \ubold{1}    & \ubold{1}     & \ubold{1}      & 0.47 & \ubold{1}     \\
  CharacterTrajectories                           & 0.85   & \ubold{1} & \ubold{1}       & \ubold{1}    & \ubold{1}  &\ubold{1}       & \ubold{1}   &\ubold{1}  \\
  Cricket                                         & \ubold{1} & \ubold{1}       & \ubold{1}    & \ubold{1}  &\ubold{1}& \ubold{1} &    \ubold{1}   &\ubold{1} \\
DuckDuckGeese                                   & 0.2    & 0.2  & 0.05    & \ubold{1}                   & 0.95             & 0.75             & -&-    \\
ERing                                           & 0.8    & 0.95 & \ubold{1}& \ubold{1}        & \ubold{1}       & \ubold{1}                  &0.85 &\ubold{1} \\
EigenWorms                                      & 0.05   & 0.3  & 0.45    & \ubold{1}                   & 0.35             & 0.52 & 0.05 &0.5\\
Epilepsy                                        & 0.05   & 0.15 & 0.3     & {0.8}                & 0.4              & 0.1        & 0.1 & \ubold{1}          \\
EthanolConcentration                            & 0.05   & 0    & 0.45    & {0.7}              & 0.05             & 0.05 & 0 &\ubold{1}                 \\
FaceDetection                                   & 0.05   & 1    & 0.65    & \ubold{0.8}           & 0.15             & 0.1          & -     &-   \\
  FingerMovements                               & \ubold{1} & \ubold{1}       & \ubold{1} & \ubold{1}                   & 0.25             & 0.05   &-&-              \\
HandMovementDirection                           & 0.05   & 0    & 0.1     & \ubold{1}      & \ubold{1}                & 0.6           &    0.1&\ubold{1}   \\
Handwriting                                 & \ubold{1} & \ubold{1}       & \ubold{1} & \ubold{1}  &  \ubold{1} & \ubold{1}            &\ubold{1}  &\ubold{1}   \\
Heartbeat                                       & 0.1    & 0    & 0.15    & \ubold{0.9}             & 0.05             & 0            & -     &-   \\
InsectWingbeat                                  & \ubold{0.2}    & 0.1  & \ubold{0.2}     & \ubold{0.2}                 & \ubold{0.2}              & 0.15   & - &-              \\
  JapaneseVowels                                 & \ubold{1} & \ubold{1}       & \ubold{1}    & \ubold{1}  &\ubold{1}& \ubold{1}& \ubold{1} &\ubold{1} \\
LSST                                            & 0      & 0.1  & 0.1     & 0.45       & 0                & 0.05         &0.14 &\ubold{1}        \\
  Libras                                          & 0.05   & 0.1  & 0.15    & 0.6                 & 0.8   & 0.55         &0.75    &\ubold{1} \\
MotorImagery                                    & 0.05   & 0    & 0.2     & \ubold{1}          & 0.7              & 0.25        & -&-\\
NATOPS                                          & 0.15   & 0.35 & 0.15    & \ubold{1}           & 0.25             & 0.2         & - &-       \\
PEMS-SF                                         & 0.4    & 0.3  & 0.4     & \ubold{1}          & 0.05             & 0.2     & - &-\\
  PenDigits                                       & \ubold{1} & \ubold{1}       & \ubold{1}    & \ubold{1}  &\ubold{1}& \ubold{1} & \ubold{1}&\ubold{1} \\
PhonemeSpectra                                  & 0      & 0.15 & 0       & 0.7         & 0.1              & 0.1      &  0.05&  \ubold{1}    \\
RacketSports                                    & 0.85   & 0.8  & 0.9     & 0.95     & 0.55             & 0.3      & 0.25   &\ubold{1}   \\
  SelfRegulationSCP1                              & 0.1    & 0.25 & 0.1     & \ubold{1}    & \ubold{1}  &\ubold{1} & 0.1 & \ubold{1}\\
  SelfRegulationSCP2                              & 0.05   & 0.1  & 0.05    & \ubold{1}    & \ubold{1}  &\ubold{1}& 0.05 &\ubold{1} \\
SpokenArabicDigits                              & 0.35   & 0.85 & 0.65    & \ubold{1}             & 0.15             & 0.2              &\ubold{1}   &\ubold{1}  \\
  StandWalkJump                                   & 0.05   & 0.6  &\ubold{1} & \ubold{1}    & \ubold{1}  &\ubold{1} &0.75 &\ubold{1}              \\
UWaveGestureLibrary                             & 0.05   & 0.05 & 0.1     & \ubold{0.15}               & 0.05             & 0          & 0.05 &0.1   \\
\end{tabular}
}
\caption{$1$ minus the Type II error for sequence length $L_0=10$ and $m=n=70$ samples.}
\label{tab:Type2L100m30}
\end{table}

\begin{table}[ht]
  \centering
  \scalebox{0.55}{
  \begin{tabular}{lcccccc}
    \toprule
    & \multicolumn{3}{c}{MMDs} & \multicolumn{3}{c}{Signature MMDs}  
    \\\cmidrule(lr){2-4}\cmidrule(lr){5-7}
    & Linear  & RBF & Laplace & Linear & RBF & Laplace\\\midrule
    ArticularyWordRecognition                       & \ubold{1} & \ubold{1}    & \ubold{1}       & \ubold{1}     & \ubold{1}         & 0.95                   \\
AtrialFibrillation                              & 0.45   & \ubold{0.9}  & 0.8     & 0.7                 & 0.7              & 0.6                  \\
    BasicMotions                                    & 0.5    & 0.9  & 0.55    & \ubold{1}      & \ubold{1}                &  \ubold{1}                    \\
    CharacterTrajectories                                & \ubold{1} & \ubold{1}    & \ubold{1}       & \ubold{1}     & \ubold{1}         & \ubold{1}               \\
    Cricket                                  & \ubold{1} & \ubold{1}    & \ubold{1}       & \ubold{1}     & \ubold{1}         & \ubold{1}                \\
DuckDuckGeese                                   & 0.15   & 0.05 & 0.2     & \ubold{0.25}                & 0.2              & 0.15                 \\
    ERing                                           & 0.85 & \ubold{1}    & \ubold{1}       & \ubold{1}     & \ubold{1}         & \ubold{1}                  \\%
EigenWorms                                      & 0.15   & 0.65 & 0.45    & \ubold{1}       & \ubold{1}        & \ubold{1}                   \\
Epilepsy                                        & 0.05   & 0.05 & 0.05    & \ubold{0.25}                & 0.1              & 0.1                  \\
EthanolConcentration                            & 0.1    & 0    & 0       & \ubold{0.15}                & 0.05             & 0.05                 \\
FaceDetection                                   & 0      & 0    & 0.15    & \ubold{0.2}                 & 0.05             & 0.1                  \\
FingerMovements                                 & \ubold{1}     & \ubold{1}    & \ubold{1}       & \ubold{1}                  & 0.95             & \ubold{1}                   \\
HandMovementDirection                           & 0.6    & 0.9  & 0.85    & \ubold{1}                  & 0.15             & 0.5                  \\
Handwriting                                     & 0.75   & \ubold{1}    & \ubold{1}       & 0.9                 & 0.95             & 0.75                 \\
Heartbeat                                       & 0.05   & 0.05 & 0       & \ubold{0.4}                 & 0.05             & 0.1                  \\
InsectWingbeat                                  & 0      & 0.05 & 0.1     & 0                   & 0.25             & 0.2                  \\
    JapaneseVowels                                  & \ubold{1}     & \ubold{1}    & \ubold{1}       & \ubold{1}                  & \ubold{1}          & \ubold{1}                   \\
LSST                                            & 0.05   & 0    & 0.05    & 0.05                & \ubold{0.1}              & 0.05                 \\
Libras                                          & 0.6    & 0.4  & 0.3     & \ubold{0.9}                 & 0.85             & 0.5                  \\
MotorImagery                                    & 0      & 0    & 0.1     & \ubold{1}          & \ubold{1}                & 0.6                  \\
NATOPS                                          & 0.3    & 0.95 & 0.95    & \ubold{1}                   & 0.2              & 0                    \\
PEMS-SF                                         & \ubold{1}      & 0.35 & \ubold{1}   & \ubold{1}        & \ubold{1}              & 0.15                 \\
    PenDigits                                     & \ubold{1} & \ubold{1}    & \ubold{1}       & \ubold{1}     & 0.95         & 0.9               \\
PhonemeSpectra                                  & 0      & 0    & 0       & \ubold{0.3}                 & 0.15             & 0                    \\
RacketSports                                    & 0.35   & 0.65 & \ubold{0.9}     & 0.5                 & 0.05             & 0.25                 \\
    SelfRegulationSCP1                                     & \ubold{1} & \ubold{1}    & \ubold{1}       & \ubold{1}     & \ubold{1}         & \ubold{1}                \\
    SelfRegulationSCP2                                     & \ubold{1} & \ubold{1}    & \ubold{1}       & \ubold{1}     & \ubold{1}         & \ubold{1}                \\
SpokenArabicDigits                              & \ubold{1}      & \ubold{1}    & \ubold{1}       & 0.95                & 0                & 0.1                  \\
StandWalkJump                                   & 0.1    & \ubold{0.95} & 0.8     & 0.65                & 0.75             & 0.6                  \\
    UWaveGestureLibrary                                    & \ubold{1} & \ubold{1}    & \ubold{1}       & 0.85     & \ubold{1}         & \ubold{1}                \\
\end{tabular}
}
\caption{$1$ minus the Type II error for sequence length $L_0=100$ and $m=n=30$ samples.}
\label{tab:Type2L200m30}
\end{table}

\begin{table}[ht]
  \centering
  \scalebox{\tablescale}{
  \begin{tabular}{lcccccc}
    \toprule
    & \multicolumn{3}{c}{MMDs} & \multicolumn{3}{c}{Signature MMDs} 
    \\\cmidrule(lr){2-4}\cmidrule(lr){5-7}
    & Linear  & RBF & Laplace & Linear & RBF & Laplace \\\midrule
ArticularyWordRecognition                       & \ubold{1} & \ubold{1}    & \ubold{1}       & \ubold{1}     & \ubold{1}         & \ubold{1}                    \\
    AtrialFibrillation                             & \ubold{1} & \ubold{1}    & \ubold{1}       & \ubold{1}     & \ubold{1}         & \ubold{1}                    \\

    BasicMotions                                   & \ubold{1} & \ubold{1}    & \ubold{1}       & \ubold{1}     & \ubold{1}         & \ubold{1}                    \\
    CharacterTrajectories    & \ubold{1} & \ubold{1}    & \ubold{1}       & \ubold{1}     & \ubold{1}         & \ubold{1}                    \\
    Cricket         & \ubold{1} & \ubold{1}    & \ubold{1}       & \ubold{1}     & \ubold{1}         & \ubold{1}                    \\
    DuckDuckGeese                                   & 0.2    & 0.2    & \ubold{1}       & \ubold{1}     & \ubold{1}         & \ubold{1}                    \\
    ERing        & \ubold{1} & \ubold{1}    & \ubold{1}       & \ubold{1}     & \ubold{1}         & \ubold{1}                    \\
    EigenWorms                                      & 0.45 & \ubold{1}    & \ubold{1}       & \ubold{1}     & \ubold{1}         & \ubold{1}                    \\
Epilepsy                                        & 0.15   & 0.75 & 0.3     &\ubold{ 0.95}                & 0.85             & 0.8                  \\
EthanolConcentration                            & 0.15   & 0.05 & 0       & \ubold{0.6}               & 0.05             & 0.1                  \\
FaceDetection                                   & 0.05   & 0.15 & \ubold{1}       & \ubold{1}   & \ubold{1}        & \ubold{1}                    \\
    FingerMovements & \ubold{1} & \ubold{1}    & \ubold{1}       & \ubold{1}     & \ubold{1}         & \ubold{1}                    \\
    HandMovementDirection    & \ubold{1} & \ubold{1}    & \ubold{1}       & \ubold{1}     & \ubold{1}         & \ubold{1}                    \\
    Handwriting      & \ubold{1} & \ubold{1}    & \ubold{1}       & \ubold{1}     & \ubold{1}         & \ubold{1}                    \\
Heartbeat                                       & 0.1    & 0.15 & 0.1     & \ubold{1}    & \ubold{1}              & \ubold{1}                  \\
InsectWingbeat                                  & 0.1    & 0.1  & 0.55    & 0.05                & \ubold{1}       & \ubold{1}                   \\
    JapaneseVowels & \ubold{1} & \ubold{1}    & \ubold{1}       & \ubold{1}     & \ubold{1}         & \ubold{1}                    \\
LSST                                            & 0.15   & \ubold{0.85} & 0.3     & 0.45                & 0.5              & 0.35                 \\
    Libras & \ubold{1} & \ubold{1}    & \ubold{1}       & \ubold{1}     & \ubold{1}         & \ubold{1}                    \\
MotorImagery                                    & 0.05   & 0.2  & 0.15    & \ubold{1}                & \ubold{1}            & \ubold{1}                    \\
    NATOPS  & \ubold{1} & \ubold{1}    & \ubold{1}       & \ubold{1}     & \ubold{1}         & \ubold{1}                    \\
    PEMS-SF  & \ubold{1} & \ubold{1}    & \ubold{1}       & \ubold{1}     & \ubold{1}         & \ubold{1}                    \\
    PenDigits & \ubold{1} & \ubold{1}    & \ubold{1}       & \ubold{1}     & \ubold{1}         & \ubold{1}                    \\
PhonemeSpectra                                  & 0.15   & 0.05 & 0.15    & \ubold{0.7}                 & 0.35             & 0.4                  \\
    RacketSports    & \ubold{1} & \ubold{1}    & \ubold{1}       & \ubold{1}     & \ubold{1}         & \ubold{1}                    \\
    SelfRegulationSCP1    & \ubold{1} & \ubold{1}    & \ubold{1}       & \ubold{1}     & \ubold{1}         & \ubold{1}                    \\
    SelfRegulationSCP2  & \ubold{1} & \ubold{1}    & \ubold{1}       & \ubold{1}     & \ubold{1}         & \ubold{1}                    \\
    SpokenArabicDigits & \ubold{1} & \ubold{1}    & \ubold{1}       & \ubold{1}     & \ubold{1}         & \ubold{1}                    \\
    StandWalkJump & \ubold{1} & \ubold{1}    & \ubold{1}       & \ubold{1}     & \ubold{1}         & \ubold{1}                    \\
    UWaveGestureLibrary  & \ubold{1} & \ubold{1}    & \ubold{1}       & \ubold{1}     & \ubold{1}         & \ubold{1}                    \\
  \end{tabular}
  }
  \caption{1 minus the Type II error for sequence length $L_0=200$ and $m=n=200$ samples.}
  \label{tab:Type2L200m200}
\end{table}

\section{Summary}

We generalized the classical moment map from the domain of random variables in finite dimensions to the domain of path-valued random variables.  
This yields a robust, universal, and characteristic feature map for stochastic processes, that can be kernelized.
The associated MMD becomes a metric for laws of stochastic processes that can be efficiently estimated from finite samples.
Let us highlight that our non-parametric approach of studying probability measures on $\bigcup_{T>0}C([0,T],\calX)$ encompasses a large class of examples which are important to different communities, e.g.
\begin{itemize}
\item genuinely discrete sequential data such as text,\footnote{If $\calX$
    is linear, sequences in $\calX$ embed into paths in $\calX$,
    $\calX^n\hookrightarrow C([0,1],\calX)$, via the Donsker embedding
  $(x_i)_{i=1}^n\mapsto \left(t\mapsto x_{\lfloor nt \rfloor}+\left( nt-\lfloor
      nt \rfloor \right)\sum_{i=1}^{\lceil nt\rceil}x_i\right)$.
  This pushes measures on sequences to measures on paths and allows to
  treat genuinely discrete data in our framework; e.g.~for text, $\calX$ is the free vector
  space spanned by letters and the path associated with a text (a sequence of
  letters) is a lattice path; multi-variate time series are sequences in $\calX=\RR^d$, etc.}
\item classical time series such as (G)ARCH, ARMA, etc.~as used in econometrics,
\item stochastic differential equations and semimartingales used in stochastic calculus,
\item stochastic processes in space and time, e.g.~stochastic partial differential equations seen as evolution equations (ODEs with noise in infinite-dimensional state space),
\item persistence paths as they arise in topological data analysis~\citep{CNO20},
\item the evolution of structured, non-Euclidean objects, such as networks, molecules, images, provided a kernel for the state space $\calX$ in which they evolve is given.
\end{itemize}
From a theoretical perspective, the existence of a universal and characteristic feature map for stochastic processes opens several research venues; two of them are 
\begin{description}

\item[Statistics in high/infinite dimensions]
  Non-parametric testing in high dimensions suffers from decreasing test power~\citep{ramdas2015decreasing}.
  Indeed, the classic kernels perform poorly in our experiments, but the robust signature kernel does well.
  This is not too surprising from a stochastic analysis perspective since sequences of iterated integrals provide a natural basis for functions of paths.
  Signatures are classical in stochastic analysis, but applications in statistics are much more recent, e.g.~\citet{papavasiliou2011parameter} use it for SDE parameter estimation and it would be interesting to connect this with the methods developed in this paper.

\item[Parametric Statistics]
  The focus of this article is the non-parametric learning of laws of stochastic processes.
  However, for some applications it can be useful to focus on a small class of stochastic processes, e.g.~certain parameterized stochastic differential equations as they appear in finance.  
   For some of these cases, semi-explicit formulas for the resulting signature moments are known (as for Brownian motion in Example~\ref{ex:BM} or more generally L{\'e}vy processes~\citep{FrizShekhar17,Chevyrev18}), or can be described by a PDE, see~\citep{ni2015expected}.
  It is natural to ask whether such strong model assumptions can be leveraged to yield signature MMDs with (semi-)explicit expressions. 
 
\end{description}
\acks{Both authors contributed equally and are listed in alphabetical order.
  IC was supported by a Junior Research Fellowship of St John's College, Oxford, when this project was started.
  HO is supported by the Oxford-Man Institute of Quantitative Finance, EPSRC grant EP/S026347/1 (Datasig), and the Alan Turing Institute.
  IC and HO would like to thank Csaba Toth for his help with the experiments and the KSig implementation of the signature MMD.}

\begin{appendices}
\section{Tensors and the Tensor Algebra\texorpdfstring{ $\TSeries{V}$}{}}\label{sec:tensor algebra}
For general background on tensors we refer to \citep[Chapter 16]{lang}.
Recall that if $V,W$ are vector spaces (possibly infinite-dimensional) then there exists another vector space $V\otimes W$ and a bilinear map $\varphi: V \times W \to V\otimes W$ with the universal property that any other bilinear map $V\times W \to Z$ factors uniquely through $\varphi$.\footnote{Concretely, if $A\subset V$ and $B\subset W$ are bases of $V$ and $W$ respectively, then $\{a\otimes b \,:\, a\in A,\;b\in B\}$ is a basis of $V\otimes W$.}
This map $\otimes$ is called the tensor product and we call the elements of $V \otimes W$ tensors.
In particular, for $m \geq 0$, we call an element of $V \otimes \cdots \otimes V=V^{\otimes m}$ a tensor of degree $m$;
by convention, $V^{\otimes 0}=\RR$.
For any vector space $V$, the space
\begin{align*}
  \TSeries{V}=   \{\bft=(\bft^m)_{m\ge 0}\,\vert\, \bft^m \in V^{\otimes m}\}
\end{align*}
is also a vector space but also carries a product. 
This product is the so-called \emph{tensor convolution product} and is defined by\footnote{Since $V^{\otimes 0} = \RR$, note that $\bfs^0,\bft^0$ are real numbers, and therefore so is $\bfs^0\bft^0\in V^{\otimes 0}=\RR$.
Furthermore, $\bfs^1\bft^0, \bfs^0\bft^1,\bfs^2\bft^0,\bfs^0\bft^2,\ldots$ are well-defined products of a tensor with a scalar.}
\begin{align}\label{eq:ncp}
  \bfs \otimes \bft \coloneqq \Big( \sum_{i=0}^m \bfs^i\otimes \bft^{m-i} \Big)_{m\geq 0} = \big( \bfs^0\bft^0, \bfs^1\bft^0 + \bfs^0\bft^1, \bfs^2\bft^0 + \bfs^1\otimes \bft^1 + \bfs^0\bft^2,\ldots \big).
\end{align}
The unit element for this product is $\bfone \coloneqq (1,0,0,\ldots)$.

The reason why $\TSeries{V}$ is a fundamental mathematical object is that $ \TSeries{V} $ is the most general algebra containing $ V $, namely $\TSeries{V}$ \emph{is the ``free algebra'' that contains $V$}, see~\citep[Chapter~1]{reutenauer-93}.
Further, it is convenient to define
\[\exp: \TSeries{V} \to \TSeries{V}, \quad \exp(\bfs) \coloneqq \sum_{m \ge 0} \frac{\bfs^{\otimes m}}{m!}.\]
(Applied with $\bfs=(0,v,0,0,\ldots)$ for $v \in \RR^d$ this reduces to the usual vector exponential).
In the main text we also use the truncated tensor algebra $\TSeriesTruncated{V}{M}$ which is again a linear space and also forms an algebra with the tensor convolution product \eqref{eq:ncp} restricted to the first $M$ tensors. 
\section{Features for Geometric Rough Paths}\label{appendix:rough}

In this appendix we give a brief introduction to rough paths theory.
For a detailed introduction see the lecture notes~\citep{lyons-04} or the monographs~\citep{friz2014course,friz-victoir-book,lyons-qian-02}.

To define the iterated integrals $\int dx^{\otimes m}$ using Riemann--Stieltjes--Young integration requires $x : [0,T] \to V$ to have finite $p$-variation for some $p\in[1,2)$.
In particular, many processes of interest in stochastic analysis fall outside this scope.
The rough paths approach is to flesh out the abstract properties of maps that associate with a path $x$ over a time interval $[s,t]$, an element of $\TSeriesTruncated{V}{\lfloor p \rfloor}$ such that it ``behaves like'' $\int
dx^{\otimes m}$.

Throughout this appendix, let $V$ be a Banach space and $p \geq 1$.
Recall from Appendix~\ref{sec:tensor algebra} that $\TSeriesTruncated{V}{M}$ is an algebra with tensor multiplication~\eqref{eq:ncp}.

\begin{definition}\label{def:RP}
A {$p$-rough path} is a continuous map
$
\mathbf{x} = (\mathbf{x}^0,\mathbf{x}^1,\ldots, \mathbf{x}^{\floor p}) : [0,T_\bfx] \to \TSeriesTruncated{V}{\lfloor p \rfloor}
$
  such that $\mathbf{x}^{0} \equiv 1$
  and, denoting\footnote{Here, for $\bft\in \TSeriesTruncated{V}{\lfloor p \rfloor}$ with $\bft^0=1$, we denote by $\bft^{-1}$ the unique element of $\TSeriesTruncated{V}{\lfloor p \rfloor}$ such that $\bft \otimes \bft^{-1}=\bfone$
(the set of all $\bft\in \TSeriesTruncated{V}{\lfloor p \rfloor}$ with $\bft^0=1$ is readily seen to be a group, so in particular $\bft^{-1}$ exists).}
 $\mathbf{x}(s,t) \coloneqq \mathbf{x}(s)^{-1}\otimes \mathbf{x}(t)$, for all $0 \leq m \leq \floor{p}$
    \begin{equation}\label{eq:finitepVar}
      \sup_{\substack{n \geq 1 \\ 0\leq t_1\leq \cdots \leq t_n \leq T_\bfx}}
      \sum_{i=1}^{n-1} \|\mathbf{x}^m(t_i,t_{i+1})\|^{p/m} < \infty.
    \end{equation}
If $x\in C([0,T_x],V)$ and there exists a $p$-rough path $\mathbf{x}$
such that 
$
\mathbf{x}^1(t) = x(t)
$,
then we call $\mathbf{x}$ a \emph{$p$-rough path lift} of $x$.
\end{definition}

For $p=1$, we recover the setting of Section~\ref{sec: ordered moments} since every $x\in \pathspVar{1}$ is a $1$-rough path by identifying $\TSeriesTruncated{V}{1} = \RR\oplus V$ with $V$.
Furthermore, for any $p \geq 1$, every $x\in \pathspVar{1}$ canonically defines a $p$-rough path $\bfx$ 
by setting $\bfx^m(t) = \int_0^t dx^{\otimes m}$, where we used the notation~\eqref{eq:mthIterInt}.
We call $\bfx$ the canonical $p$-rough path lift of $x$.

A $p$-rough path $\mathbf{x}$ is called \emph{geometric} if there exist $p$-rough paths $(\mathbf{x}_n)_{n \geq 1}$ which are the canonical lifts of smooth paths such that $\mathbf{x}_n \rightarrow \mathbf{x}$ in the $p$-variation metric $d_{p\Var}$.\footnote{See~\citep[Eq.~(3.70)]{lyons-qian-02} for the definition of the $p$-variation metric; for brevity, we do not give a definition here, particularly because its explicit form plays no role in the sequel.}

\begin{example}\label{Ex:rough path lifts}
  \begin{enumerate}[label=(\arabic*)]
\item \label{point:semimart}\textbf{(Semi)martingales.}
  With $2<p< 3$ we cover the theory of continuous (semi)martingales~\citep[Sec.~14]{friz-victoir-book}: every continuous semimartingale $X : [0,T] \to \RR^d$ is of bounded $p$-variation
  for any $p>2$.
  Stochastic  (It{\^o} or Stratonovich) integration can be used to the define the
  first $2=\lfloor p \rfloor$ iterated integrals
$
    \mathbf{X}(t)= (1,\int_0^t dX,\int_0^t dX^{\otimes 2})
$.
One can verify that~\eqref{eq:finitepVar} holds almost surely, hence
$\mathbf{X}$ is a $p$-rough path.
Moreover, if Stratonovich integration is used, then $\mathbf{X}$ is geometric.
Extensions exist to the discontinuous setting~\citep{CF19}.

\item\textbf{Gaussian processes.} Many Gaussian processes admit canonical lifts to geometric $p$-rough paths.
A criterion for such a lift to exist can be given in terms of the covariance function~\citep{coutin-qian-02, FGGR16}, which in particular covers fractional Brownian motion with Hurst parameter $H>1/4$.

\item \textbf{Markov processes.}
  Likewise, many Markov processes admit canonical lifts to geometric $p$-rough paths (typically $2 < p < 3$).
These include diffusions on fractals~\citep{bass-hambly-lyons-02} and Markov processes arising from elliptic Dirichlet forms~\citep{friz-victoir-subelliptic, CO18}.
Such processes are typically not semimartingales and thus fall outside the scope of It{\^o}--Stratonovich calculus.
\end{enumerate}
\end{example}


To mimic the construction in Section~\ref{sec: ordered moments}, we define the space of geometric $p$-rough paths
\begin{equation*}
  \pathspVar{p} \coloneqq \bigcup_{T>0}
\Big\{
\mathbf{x} : [0,T] \to \TSeriesTruncated{V}{\lfloor p \rfloor} : \mathbf{x} \text{ is a geometric $p$-rough path}
\Big\}\;,
\end{equation*}
\label{page ref pRPs}
equipped with the (non-Hausdorff) topology induced by
\begin{equation}\label{eq:rho_p}
\rho_p(\mathbf{x},\mathbf{y}) = \inf_{\tau} d_{p\Var}(\mathbf{x}, \mathbf{y}\circ\tau)\;, 
\end{equation}
where the infimum is over all increasing bijections $\tau : [0,T_\mathbf{x}] \to [0,T_{\mathbf{y}}]$.

\begin{remark}
In the case $V= \RR^d$, a geometric $p$-rough path satisfies the shuffle identity~\eqref{eq:shuffle}.
\end{remark}

\subsection{Tree-like Equivalence}\label{subsec:trees}
For a topological space $\calX$, a function $x : [0,T_x] \to \calX$ is called \emph{tree-like} if $x$ is continuous and there exists an $\RR$-tree $\frakT$, a continuous function $\phi : [0,T_x] \to \frakT$, and a map $\psi : \frakT \to \calX$ such that $\phi(0)=\phi(T_x)$ and $x = \psi\circ\phi$.
Let $\overleftarrow{x} : [0,T_x] \to \calX, \overleftarrow{x}(t) = x(T_x-t)$ denote the time-reversal of $x$.
For another function $y : [0,T_y] \to \calX$, we denote the concatenation of $x$ with $y$ by
\[
x*y : [0,T_x+T_y] \to \calX\;, \quad x*y(t) =
\begin{cases}
x(t) &\text{ if } t\in[0,T_x]
\\
y(t) &\text{ if } t \in (T_x,T_x+T_y]\;.
\end{cases}
\]
Note that $x*y$ is a continuous path if $x$ and $y$ are continuous and $x(T_x)=y(0)$.
We say that $x$ and $y$ are \emph{tree-like equivalent}, and write $x\sim_t y$, if $x*\overleftarrow y$ is tree-like.\label{page ref treelike}

\begin{example}
An important example of tree-like equivalence is when $x$ is a reparameterization of $y$ i.e.~there exists an increasing bijection $\tau : [0,T_x]\to [0,T_y]$ such that $x=y\circ\tau$.
\end{example}

\begin{example}
Another way in which paths can be tree-like equivalent is if they differ only through back-tracking excursions.
For example the following two paths $x,y\colon[0,4]\to \RR^2$ depicted in the figures below are tree-like equivalent:
$x(t)\coloneqq[t,t]$ and
\[
y(t) \coloneqq
\begin{cases}
[2t,2t] & \textnormal{ if } t\in [0,1)\\
[2+2(t-1), 2-2(t-1)] & \text{ if } t\in [1,2)\\
[4-2(t-2), 2(t-2)] & \text{ if } t\in [2,3)\\
[2+2(t-3), 2 + 2(t-3)] & \text{ if } t\in [3,4]\;.\\
\end{cases}
\]
\[
      \begin{tikzpicture}
          \begin{axis}[
			name=Xtree,
			  unit vector ratio*=1 1 1,
			  disabledatascaling,
              height       = 1.9in,
              xmax         = 4.3,
xmin         = -0.3,
              ymax         = 4.3,
ymin         = -0.3,
			  grid=both,
              xtick        = {-1,0.0, 1,2, 3.0,4},
              xticklabels  = {$-1$, $0$, $1$,$2$,$3$,$4$},
              ytick        = {0.0, 1,2,3,4 },
              yticklabels  = {0, $1$,$2$,$3$,$4$},
              axis lines   = center,
              line cap=round,
              ]
\addplot[->,color=red,very thick] coordinates { (0,0) (4,4) };
          \end{axis}
\node[anchor=north] at (-0.7,2) {$x:$};
      \end{tikzpicture}
\qquad\qquad
\begin{tikzpicture}
          \begin{axis}[
name=Ytree,
			  unit vector ratio*=1 1 1,
			  disabledatascaling,
              height       = 1.9in,
              xmax         = 4.3,
xmin         = -0.3,
              ymax         = 4.3,
ymin         = -0.3,
			  grid=both,
              xtick        = {-1,0.0, 1,2, 3.0,4},
              xticklabels  = {$-1$, $0$, $1$,$2$,$3$,$4$},
              ytick        = {0.0, 1,2,3,4 },
              yticklabels  = {0, $1$,$2$,$3$,$4$},
              axis lines   = center,
              line cap=round,
              ]
\addplot[-,color=red,very thick] coordinates { (0,0) (2,2) };
\addplot[->,color=red,very thick] coordinates { (2,2) (4,0)};
\addplot[->,color=red,very thick] coordinates { (4,0) (2,2)};
\addplot[->,color=red,very thick] coordinates { (2,2) (4,4)};
          \end{axis}
\node[anchor=north] at (-0.7,2) {$y:$};
      \end{tikzpicture}
\]
\end{example}
 
\subsection{Ordered Moments for Rough Paths}\label{subsec:liftingRP}
A key result of~\citet{lyons-98} is that the signature map $\Sig$ defined by~\eqref{eq:ordered moments map BV} extends to the space of $p$-rough paths: for any $\bfx\in
\pathspVar{p}$, the integrals $\int d\bfx^{\otimes m}$ are canonically defined for all $m > \floor{p}$.
Moreover, $\Sig$ is injective on $\pathspVar{p}$ up to tree-like equivalence~\citep{BGLY16}.
We summarise these results in the following extension of Theorem~\ref{thm:propsPhiMap}.
\begin{theorem}\label{Thm:extension}
There exists a map
$\Sig:\pathspVar{p}\rightarrow \FSeriesOne{V}$, $\Sig(\bfx) = (1, \Sig^1(\bfx),\Sig^2(\bfx), \ldots) $
such that $\Sig^m(\mathbf{x}) = \mathbf{x}^m(T_\bfx)$ for all $m=0,\ldots, \floor{p}$.
Furthermore,
$\Sig(\mathbf{x}) = \Sig(\mathbf{y})$ if and only if $\bfx\sim_t\bfy$.
\end{theorem}

\begin{example}\quad
  \begin{itemize}
  \item 
  \label{ex:BVasRP}
  Let $x \in \pathspVar{1}$.
  Then $\Sig(x) \in \FSeriesOne{V}$ defined by~\eqref{eq:ordered moments map BV} agrees with $\Sig(x)$ in Theorem~\ref{Thm:extension}.
\item
Let $X=(X(t))_{t\in[0,T]}$ be a continuous semimartingale  in $\RR^d$ with geometric $p$-rough path lift
$\mathbf{X}=\left( 1,\mathbf{X}^1,\mathbf{X}^2 \right)\equiv \left( 1,\int dX,
  \int dX^{\otimes 2} \right)$ as in part~\ref{point:semimart} of
Example~\ref{Ex:rough path lifts}.
Then
\[
\Sig(\mathbf{X}) = \Big(1, \int_0^T dX, \int_0^T dX^{\otimes 2},\int_0^T
dX^{\otimes 3}, \ldots \Big)\in\FSeriesOne{\RR^d},
\]
where the stochastic integrals are taken in the Stratonovich sense.
\end{itemize}
\end{example}

Theorem~\ref{Thm:extension} suggests the following generalization of Definition~\ref{def:BV tree-like equivalence}.
\begin{definition}\label{def:RP tree-like equivalence}
We define the space of unparameterized geometric $p$-rough paths $\pathsreduced{p}$ as the set of equivalence classes $\pathspVar{p} / \sim_t$.\label{page ref RPs unpar} We equip $\pathsreduced{p}$ with the topology induced by the signature map embedding $\Sig : \pathsreduced{p} \to \FSeriesOne{V}$.
\end{definition}


\subsection{Time Parameterization}
\label{subsec:time par RPs}

As in the bounded variation setting of Section~\ref{subsec:timepar BV}, we can make the map $\Sig$ variant to the parameterization of time by adding a time component.
For any $\bfx\in\pathspVar{p}(V)$, there is a canonical $p$-rough path lift of $t \mapsto (\bfx^1(t),t)$, denoted by $\overline \bfx\in\pathspVar{p}(V\oplus \RR)$, and which extends $\bfx$ in the sense that $\langle \ell,\bfx^m\rangle = \langle\ell, \overline\bfx^m\rangle$ for all $m=0,\ldots, \floor p$ and $\ell \in (V^{\otimes m})'$, see~\citep[Sec.~3.3.3]{lyons-qian-02}.
As such, we can identify $\pathspVar{p}(V)$ with a subset of $\pathspVar{p}(V\oplus \RR)$.
Furthermore, $\overline \bfx \sim_t \overline \bfy$ if and only if $\bfx=\bfy$, and thus $(\pathspVar{p}(V),\overline\rho)$ is a genuine metric space, where $\overline\rho_p$ is the quantity~\eqref{eq:rho_p} associated to $\pathspVar{p}(V\oplus \RR)$.

\subsection{Topological Properties and MMD for Measures on Rough Paths}
\label{appendix:MMD for RPs}

\begin{lemma}\label{lem:complete}
Let $(\bfx_n)_{n=1}^\infty$ be a Cauchy sequence for $\rho_{p}$. Then $\rho_p(\bfx_n,\bfx) \to0$ for some $\bfx \in \pathspVar{p}$.
\end{lemma}

\begin{proof}
Denoting by $\pathspVar{p}_{[0,T]}$ the space of geometric rough paths defined on $[0,T]$, recall that $(\pathspVar{p}_{[0,T]},d_{p\Var})$ is a complete metric space~\citep[Lem.~3.3.3]{lyons-qian-02}.
We can find a subsequence $n(k)$ 
and time-changes $\tau_k : [0,1]\to [0,T_{\bfx_{n(k)}}]$ such that, denoting $\bfy_k = \bfx_{n(k)}\circ \tau_k\in \pathspVar{p}_{[0,1]}$, $d_{p\Var}(\bfy_k,\bfy_{k+1}) < 2^{-k}$.
In particular, $\bfy_k$ is Cauchy for $d_{p\Var}$ and converges to some $\bfx\in \pathspVar{p}_{[0,1]}$, which implies that $\rho_p(\bfx_n,\bfx) \to 0$. 
\end{proof}

\begin{proposition}\label{prop:Polish}
The topological space $\pathsreduced{p}$ is metrizable and the metric space $(\pathspVar{p},\overline\rho_p)$ is complete.
Moreover, if $V$ is separable, then so are $(\pathspVar{p},\rho_p)$ and $(\pathspVar{p},\overline\rho_p)$.
\end{proposition}

\begin{proof}
Since $\Sig$ is injective and $\FSeriesOne{V}$ is a metric space, $\pathsreduced{p}$ is metrizable.
Completeness of $(\pathspVar{p},\overline\rho_p)$ follows from Lemma~\ref{lem:complete} and the fact that
$\{\overline\bfx : \bfx\in \pathspVar{p}\}$ is a closed subset of $(\pathspVar{p}(V\oplus\RR,\overline\rho_p)$.
If $V$ is separable, then $C^\infty([0,T],V)$ is separable under the $C^1$-norm $\|\cdot\|_{C^1}$, see~\citep[Sol.~2.13]{friz2014course}.
Since the canonical lift map $(C^\infty,\|\cdot\|_{C^1})\ni x\mapsto \bfx\in(\pathspVar{p},\rho_{p})$ is continuous, separability of $(\pathspVar{p},\rho_p)$ and $(\pathspVar{p},\overline\rho_p)$ follows.
\end{proof}

We conclude this section with the proof of Theorem~\ref{thm: signature kernel}\ref{point:topology comparison}.
Let notation be as in Section~\ref{subsec:kernel mean embedding} and Theorem~\ref{thm: signature kernel}. 
\begin{proposition}\label{prop:weaker than weak}
\begin{enumerate}[label=(\roman*)]
  \item \label{point:not weak conv} There exist probability measures $\mu_n,\mu$, on $\pathsreduced{p}$ such that $d_{\kernel}(\mu_n,\mu)\to 0$ but such that $\mu_n$ does not converge weakly to $\mu$.
\item \label{point:strictly weaker} Suppose $H$ is finite dimensional
  and that $\calM$ is a set of probability measures on $\pathsreduced{p}$ which is compact under weak convergence.
  Then $d_\kernel$ and weak convergence induce the same topology on $\calM$.
  In particular, weak convergence implies convergence in $d_{\kernel}$.
\end{enumerate}
The same statements hold with $\pathsreduced{p}$ replaced by $(\pathspVar{p},\overline\rho_p)$.
\end{proposition}
\begin{proof}
\ref{point:not weak conv} Let $p < p' < \floor{p}+1$.
Consider $\mathbf{x} \in \pathsreduced{p}$ and a sequence $\mathbf{x}_n \in \pathsreduced{p}$ such that $\mathbf{x}_n$ does not converge to $\mathbf{x}$ as elements of of $\pathsreduced{p}$ but such that $\mathbf{x}_n\to \mathbf{x}$ as elements of $\pathsreduced{p^\prime}$, i.e., in the $p^\prime$-variation metric (it is a simple exercise to construct such a sequence for any $\mathbf{x} \in \pathsreduced{p}$).
Due to continuity of $\Phi : \pathsreduced{p^\prime} \to \FSeries{H}$, we conclude that $\|\Phi(\mathbf{x}_n) - \Phi(\mathbf{x})\| \to 0$.
In particular, for the corresponding Dirac delta measures, we have
\[
d_{\kernel}(\delta_{\mathbf{x}_n}, \delta_{\mathbf{x}})
= \sup_{\mathbf{t} \in \FSeries{H}, \|\mathbf{t}\| \leq 1} |\langle\mathbf{t},\Phi(\mathbf{x}_n)\rangle -|\langle\mathbf{t},\Phi(\mathbf{x})\rangle |
\to 0.
\]
However, due to the assumption that $\mathbf{x}_n$ does not converge to $\mathbf{x}$ in $\pathsreduced{p}$, it holds that $\delta_{\mathbf{x}_n}$ does not converge weakly to $\delta_{\mathbf{x}}$ as probability measures on $\pathsreduced{p}$.

\ref{point:strictly weaker}
Let $\mu$ and $\mu_n$ be probability measures on $\pathsreduced{p}$ such that $\mu_n \to \mu$ weakly.
Observe that
\begin{align*}
d_{\kernel}(\mu_n,\mu)^2
  &= \sup_{\mathbf{t} \in \FSeries{H}, \|\mathbf{t}\| \leq 1} \Big| \int_{\pathsreduced{p}} \langle \mathbf{t}, \Phi(\mathbf{x}) \rangle \mu_n(d\mathbf{x}) - \int_{\pathsreduced{p}} \langle \mathbf{t}, \Phi(\mathbf{x}) \rangle \mu(d\mathbf{x})\Big|^2
\\
  &= \Big\|\int_{\pathsreduced{p}}\Phi(\mathbf{x})\mu_n(d\mathbf{x})	- \int_{\pathsreduced{p}}\Phi(\mathbf{x})\mu(d\mathbf{x}) \Big\|^2_{\FSeries{H}}
\\
  &= \sum_{m \geq 0} \Big\| \int_{\pathsreduced{p}} \lambda(\Sig(\bfx))^m \Sig^m(\mathbf{x}) \mu_n(\mathbf{x}) - \int_{\pathsreduced{p}}\lambda(\Sig(\bfx))^m \Sig^m(\mathbf{x}) \mu(\mathbf{x}) \Big\|^2_{H^{\otimes m}}.
\end{align*}
Note that for every $m \geq 0$, $\mathbf{x} \mapsto \lambda(\Sig(\bfx))^m \Sig^m(\mathbf{x})$ is a continuous bounded function from $\pathsreduced{p}$ into the finite dimensional vector space $H^{\otimes m}$.
It follows from the weak convergence $\mu_n \to \mu$ that the final quantity converges to zero as $n \to \infty$, thus weak convergence implies convergence in $d_\kernel$.
Therefore $\calM$ is also compact under $d_\kernel$,
and since the topologies induced by $d_\kernel$ and weak convergence are comparable, the two topologies necessarily coincide on $\calM$.
Therefore $\calM$ is also compact under $d_\kernel$,
and since the topologies induced by $d_\kernel$ and weak convergence are comparable, the two topologies necessarily coincide on $\calM$.

The corresponding claims for $\pathspVar{p}$ follow in an identical manner.
\end{proof}
\subsection{Example: Non-robust Signature Moments gone wrong}\label{ex:example nonchar}
Define $X,Y\colon[0,1]\to\RR^2$ by $X_t=t \cdot N^\top $ and $Y_t= t \cdot M^\top$, where $N=(N_1,N_2)$ consists of two-independent lognormal distributions and $M=(M_1,M_2)$ consists of two independent perturbed lognormals, i.e.~$p(n_1,n_2) = \prod_{i=1}^n\frac{1}{n_i\sqrt{2\pi}} \exp (\frac{-\log^2 (n_i)}{2})$ and $q(m_1,m_2)=p(m_1,m_2) \prod_{i=1}^2 (1+\sin (2 \pi \log m_i) )$. 
Then all signature moments of $X$ and $Y$ coincide although $X$ and $Y$ have different laws,
  \begin{align}\label{eq:same moments}
    \EE[\int dX^{\otimes m} ]  = \EE[ \int dY^{\otimes m}] \text{ for all }m\ge0\;.
  \end{align}
  The moment equality \eqref{eq:same moments} holds since a direct calculation shows that $\int dX^{\otimes m}= \frac{(X(1)-X(0))^{\otimes m}}{m!}=\frac{N^{\otimes m}}{m!}$ and analogously $\int dY^{\otimes m} =\frac{M^{\otimes m}}{m!}$.
  Hence the only information captured in $\EE[\int dX^{\otimes m}]$ and $\EE[\int dY^{\otimes m}]$ are the $m$-th moments of $M$ and $N$, but $\EE[M^{\otimes{ m}}] = \EE[N^{\otimes{ m}}] $ follows from a direct calculation. 
\section{Kernel Background}
\label{appendix:kernels}
Throughout this appendix, we fix a topological space $\calX$, an inner product space $(E,\langle \cdot,\cdot \rangle_E)$ and feature map $\Phi : \calX \to E$.
For $x \in \calX$, we define the kernel function $k_x \in \RR^{\calX}$ by
$k_x(y) \coloneqq k(x,y)$, where $k(x,y) \coloneqq\langle \Phi(x),\Phi(y) \rangle$.
Consider the subspace
$
\calH_0 \coloneqq \operatorname{span}\{k_x:x\in \calX\} \subset \RR^\calX
$
equipped with the inner product defined by $\langle k_x,k_y \rangle_{\calH_0} \coloneqq k(x,y)$.
We first recall the following theorem of Moore--Aronszajn.
\begin{theorem}[Moore--Aronszajn]
There exits a unique Hilbert space $\calH\subset \RR^\calX$ with $k$ as the reproducing kernel.
Moreover, $\calH_0$ is dense in $\calH$.
\end{theorem}

The following theorem clarifies the relation to our original feature map $x\mapsto\Phi(x)$ and $x\mapsto k_x$ and shows how to construct $\calH$ as a subspace of (the completion of) $E$.
\begin{theorem}\label{Thm:features2kernel}
\begin{enumerate}[label=(\roman*)]
\item\label{item:HdenseE} There exists a unique linear map $\Psi : \calH_0 \to E$ such that the following diagram commutes:
\[
  \begin{tikzcd}
    \calX \arrow{r}{x \to k_x} \arrow[swap]{dr}{\Phi} & \calH_0 \arrow{d}{\Psi} \\
     & E\;.
  \end{tikzcd}
\]
\item\label{item:PsiInj} The map $\Psi$ is injective and is an isometry onto its image.
\item\label{item:iPhi}
Denote by $x \mapsto x^\prime$ the canonical map $E \to E^\prime$ which identifies $E$
with (a subspace of) $E^\prime$.
Then
\[
  \imath(\Phi(x)^\prime) = k_x
\]
and the following diagram commutes:
     \[
  \begin{tikzcd}
    \calX \arrow[r, "x \mapsto k_x"] \arrow[swap, dr, "\Phi"] & \calH_0 \arrow[d, "\Psi"] \arrow[r, hook, "\operatorname{id}"] & \RR^\calX  \\
     & E \arrow[r, "x \mapsto x'"] & E' \arrow[u, "\imath"]\;.
  \end{tikzcd}
\]
\item \label{item:dense} 
The image of $\calH_0$ in $E^\prime$ under $h \mapsto \Psi(h)^\prime$ is a dense subspace of $\Ker(\imath)^\perp$.
\end{enumerate}
\end{theorem}

\begin{proof}
For Point~\ref{item:HdenseE}, existence and uniqueness of $\Psi$ follows from the observation that if $\sum a_i k_{x_i} \equiv 0$, then $\sum a_i \langle \Phi(x_i), \Phi(y) \rangle = 0$ for all $y \in \calX$, and thus $\sum a_i \Phi(x_i)$ is an element of both $\Phi(\calX)^\perp$ and $\operatorname{span}\left[\Phi(\calX)\right]$, and thus must be zero.
For Point~\ref{item:PsiInj}, note that $\Psi(\sum a_i k_{x_i}) = 0$ is equivalent to $\sum a_i \Phi(x_i) = 0$, so that
\[
\sum a_i \langle \Phi(x_i), \Phi(y)\rangle = \sum a_i k_{x_i}(y) = 0, \; \; \forall y \in \calX.
\]
It follows that $\Psi$ is injective. The fact that $\Psi$ is an isometry follows from Point~\ref{item:HdenseE}.
Point~\ref{item:iPhi} follows immediately since for all $x,y \in \calX$
\[
\imath(\Phi(x)')(y) = \Phi(x)' \circ \Phi(y) = \langle \Phi(x), \Phi(y) \rangle = k(x,y) = k_x(y).
\]
For Point~\ref{item:dense}, we define for any subset $F \subset E$ the set $F^\circ \coloneqq \{f' \in E' \mid f'(f) = 0, \forall f \in F\}$. 
Now $\Ker(\imath)$ consists of all $z \in E^\prime$ such that $z^\prime(\Phi(x)) = 0$ for all $x \in \calX$, so that $\Ker(\imath) = \Phi(\calX)^\circ \subset E^\prime$.
Since $\operatorname{span}\left[\Phi(\calX)\right] = \Psi(\calH_0)$, it follows that $\Ker(\imath) = \Psi(\calH_0)^\circ$. The conclusion now follows from the fact that for any subspace $F \subset E$, the image of $F$ under $x \mapsto x'$ is dense in $(F^\circ)^\perp$.
\end{proof}
\begin{proposition}\label{prop:kPhiUnivCharFull}
Suppose $\calF \subset \RR^\calX$ is a locally convex TVS and that $\calH_0$ continuously embeds into $\calF$.
Then the map $\imath : E' \to \RR^\calX$ given by~\eqref{eq:iEmbedding} maps $E'$ continuously into $\calF$.
Furthermore 
        \begin{itemize}
        \item $\Phi$ is universal to $\calF$ iff the kernel $k$ is universal to $\calF$,
        \item $\Phi$ is characteristic to $\calF^\prime$ iff the kernel $k$ is characteristic to $\calF^\prime$.
       \end{itemize}
\end{proposition}
\begin{proof}
Substituting $\calF$ by its completion if necessary, we may assume $\calF$ is complete.
Write $F_0 \coloneqq \Psi(\calH_0)' \subset E'$ and let $F \subset E'$ denote the closure of $F_0$ in $E'$.
By Point~\ref{item:iPhi} of Theorem~\ref{Thm:features2kernel}, it holds that
\begin{equation}\label{eq:imagesEqual}
\id(\calH_0) = \imath(F_0) \subset \calF,
\end{equation}
so by the assumption that $\id : \calH_0 \hookrightarrow \calF$ is continuous, the restriction $\imath|_{F_0} : F_0 \to \calF$ is continuous.
By definition of $\imath$, it is easy to see that the restriction $\imath|_{F}$ agrees with the unique continuous extension of $\imath|_{F_0}$ to $F$.
Hence $\imath\mid_F : F \to \calF$ is continuous.
We now write $E' = F \oplus F^\perp$.
By Point~\ref{item:dense} of Theorem~\ref{Thm:features2kernel}, we have $F = \Ker(\imath)^\perp$.

Note that $\Ker(\imath)$ is closed in the weak topology of $E^\prime$, and thus, \textit{a fortiori}, under the strong (norm) topology.
Indeed, if $\ell_n \rightarrow \ell$ pointwise in $E^\prime$, then for all $ x\in \calX$, $\lim_{n \rightarrow \infty}\ell_n(\Phi(x)) = \ell(\Phi(x))$. In particular, if $\ell_n \in \Ker(\imath)$ for all $n \geq 1$, then $\ell \in \Ker(\imath)$.

As a consequence, it holds that $\Ker(\imath) = F^\perp$, and thus
\begin{equation}\label{eq:EdecompF}
E' = F \oplus \Ker(\imath),
\end{equation}
from which the continuity of $\imath : E' \to \calF$ follows.

Finally, by~\eqref{eq:EdecompF} and the continuity of $\imath$, the closures in $\calF$ of $\imath (F_0)$ and $\imath(E')$ coincide. Combining with~\eqref{eq:imagesEqual}, it holds that $\imath(E')$ is dense in $\calF$ iff $\id(\calH_0)$ is dense in $\calF$, from which the first equivalence follows. The second equivalence now follows from Theorem~\ref{thm:univIffChar} and Proposition \ref{prop:kernel mean iff char}.
\end{proof}

\subsection{Signature Kernel Discretization}
\label{appendix:sig_kernels}

Suppose we are in the setting of Section~\ref{subsec:comp_sig_kernel}.

\begin{proposition}\label{prop:kernel_discrete}
\begin{enumerate}[label=(\roman*)]
\item \label{point:seq_kernel} The kernel $ \kernelTrunc{M}^{\kappa,+}:  \calX^+\times  \calX^+\rightarrow \RR$ is a bounded, positive semidefinite kernel for $\calX^+$.

Moreover, for $x,y \in \pathspVar{1}(\calX)$ and partitions $\pi \subset [0,T_x],\pi'\subset[0,T_y]$
\begin{equation}\label{eq:discrete_approx}
\big| \kernel^{\kappa}(x,y)- \kernelTrunc{M}^{\kappa,+}(x^\pi, y^{\pi'}) \big|
\leq
\overline K
(F(Q_x)+F(Q_y) 
+
F(\Delta^\pi_{x})
+
F(\Delta^{\pi'}_y))\;,
\end{equation}
where $x^\pi$ (resp. $y^{\pi'}$) denote the sequences given by sampling $x$, $y$ along $\pi$ (resp.~$\pi'$) and
\begin{align*}
  \Delta^\pi_{x}&\coloneqq \|\kappa_x\|_{1\Var}e^{\|\kappa_x\|_{1\Var}}\max_{t_i\in \pi}\|\kappa_x\|_{1\Var;[t_{i-1},t_{i}]}\;,\\
Q_x &\coloneqq \frac{e^{\|\kappa_x\|_{1\Var}}\|\kappa_x\|_{1\Var}^{M+1}}{(M+1)!}\\
\overline K &\coloneqq \|\psi\|_\infty^{1/2}(1+K^{1/2}+\sqrt 2\|\psi\|_\infty^{1/2})\;,\\
F(x)&\coloneqq x\vee\sqrt x\;.
\end{align*}
\item \label{itm:sequence MMD}
  Consider probability measures $\mu$,$\nu$ on $\pathspVar{1}(\calX)$ whose marginals on $\calX^+$ along partitions $\pi$,$\pi'$ are $\mu^\pi$,$\nu^{\pi'}$. 
  Then
  \[
    |d_{\kernel^\kappa}(\mu,\nu)-d_{\kernelTrunc{M}^{\kappa,+}}(\mu^\pi,\nu^{\pi'})| \leq 4\overline K
    \EE[F(Q_X) + F(\Delta^\pi_X) + F(Q_Y) + F(\Delta^{\pi'}_Y)]\;.
  \]
The same inequality applies if the partitions are also random (see the proof for the formal statement). 

\end{enumerate}
\end{proposition}

\begin{proof}
  For brevity we denote 
   $\Sig_{M}^{\linear}(h) \coloneqq \prod_{i=1}^\ell (1+h_i - h_{i-1})$.
 \citep[Cor.~4.3]{OK16} implies that
 \begin{align}\label{eq:linear signature}
   \|\Sig^{\linear}_{M}(\kappa_{x^\pi}) - \Sig_M(\kappa_x)\|
 \leq \Delta^\pi_{x}.
\end{align}

The first part of \ref{point:seq_kernel} simply follows from the boundedness of $\Lambda$ and the fact that an inner product, in this case on $\TSeriesTruncated{H}{M}$, is positive definite. 
The remaining estimate~\eqref{eq:discrete_approx} follows by using \eqref{eq:linear signature} in the technical Lemmas~\ref{lem:CS} and~\ref{lem:untrunc_to_trunc} that quantifies the factorial decay in the truncation level $M$.

For point \ref{itm:sequence MMD} we directly prove the more general statement that allows for random partitions.
Therefore denote $\Pi \coloneqq \cup_{T>0} \Pi(T)$ with $\Pi(T)\coloneqq\{\pi:\pi=\{0 \le t_1 < \cdots \le T\}\}$, and set $\calS \coloneqq \{(x,\pi)\in \pathspVar{1}(\calX)\times \Pi : \pi\in\Pi(T_x)\}$.
We define a kernel on $\calS$ as $\kernelTrunc{M}^{\kappa,\Pi}((x,\pi),(y,\pi')) \coloneqq \kernelTrunc{M}^{\kappa,+}(x^\pi,y^{\pi'})$.
Note that the associated MMD $d_{\kernelTrunc{M}^{\kappa,\Pi}}$ equals $d_{\kernelTrunc{M}^{\kappa,+}}$ for the special case of deterministic partitions.
The claim then follows by using the representation of MMD as a sum of expectations of kernels (equation \eqref{eq:MMD_expectations}), and point~\ref{point:seq_kernel}.
\end{proof}

In the proof of Proposition~\ref{prop:kernel_discrete}, we used the following two lemmas. We keep notation as in Proposition~\ref{prop:kernel_discrete}.

\begin{lemma}\label{lem:CS}
For any $\bfs,\bar \bfs,\bft,\bar \bft\in \FSeriesOne{H}$
\begin{equation}
|\scal{\Lambda(\bfs),\Lambda(\bft)} - \scal{\Lambda(\bar \bfs),\Lambda(\bar \bft)}|
\leq \overline K(
F(\|\bfs - \bar \bfs\|)+
F(\|\bft - \bar \bft\|)
)\;.
\end{equation}
\end{lemma}

\begin{proof}
This follows from Proposition~\ref{prop:LipschitzBoundPhi}\ref{point:maps into ball}-\ref{point:Lipschitz}
and the Cauchy--Schwarz inequality.
\end{proof}

Setting $\varphi(x)=\kappa(x,\cdot)$ as in Proposition~\ref{prop:lift_kernel_univ},
recall the space $\pathspVar{1}(\calX)$ from Definition~\ref{def:normSigLift}.

\begin{lemma}\label{lem:untrunc_to_trunc}
Let
$x,y \in \pathspVar{1}(\calX)$.
Then
\begin{equation}\label{eq:trunc_approx}
\big| \kernel^\kappa(x,y) - \kernelTrunc{M}^{\kappa}(x,y) \big|
\leq
\overline K
(F(Q_x)+F(Q_y) )\;.
\end{equation}
\end{lemma}

\begin{proof}
  Denoting $\bft = (\bft^0,\bft^1,\ldots) \coloneqq \Sig(\kappa_x) \in\FSeriesOne{H}$,
recall that $\|\bft^m\| \leq \|\kappa_x\|_{1\Var}^m/m!$~\citep[Prop.~2.2]{lyons-04}.
Let $\pi_M \colon \FSeries{H} \to \TSeriesTruncated{H}{M}$ denote the level-$M$ truncation map.
Then
\[
\|\pi_M \bft - \bft\| \leq
\Big(\sum_{m>M} \Big(\frac{\|\kappa_x\|_{1\Var}^{m}}{m!}\Big)^2\Big)^{1/2}
\leq
\sum_{m>M} \frac{\|\kappa_x\|_{1\Var}^{m}}{m!} \leq
\frac{e^{\|\kappa_x\|_{1\Var}} \|\kappa_x\|_{1\Var}^{M+1}}{(M+1)!} = Q_x\;.
\]
The conclusion now follows from Lemma~\ref{lem:CS}.
\end{proof}

\end{appendices}

\bibliography{roughpaths}

\def\cprime{$'$}
\begin{thebibliography}{53}
\providecommand{\natexlab}[1]{#1}
\providecommand{\url}[1]{\texttt{#1}}
\expandafter\ifx\csname urlstyle\endcsname\relax
  \providecommand{\doi}[1]{doi: #1}\else
  \providecommand{\doi}{doi: \begingroup \urlstyle{rm}\Url}\fi

\bibitem[Bagnall et~al.(2018)Bagnall, Dau, Lines, Flynn, Large, Bostrom,
  Southam, and Keogh]{bagnall2018uea}
Anthony Bagnall, Hoang~Anh Dau, Jason Lines, Michael Flynn, James Large, Aaron
  Bostrom, Paul Southam, and Eamonn Keogh.
\newblock The {UEA} multivariate time series classification archive, 2018.

\bibitem[Bass et~al.(2002)Bass, Hambly, and Lyons]{bass-hambly-lyons-02}
Richard Bass, Ben Hambly, and Terry Lyons.
\newblock Extending the {W}ong-{Z}akai theorem to reversible {M}arkov
  processes.
\newblock \emph{J. Eur. Math. Soc. (JEMS)}, 4\penalty0 (3):\penalty0 237--269,
  2002.

\bibitem[Berlinet and Thomas-Agnan(2011)]{berlinet2011reproducing}
Alain Berlinet and Christine Thomas-Agnan.
\newblock \emph{Reproducing kernel Hilbert spaces in probability and
  statistics}.
\newblock Springer Science \& Business Media, 2011.

\bibitem[Berndt and Clifford(1994)]{berndt1994using}
Donald~J Berndt and James Clifford.
\newblock Using dynamic time warping to find patterns in time series.
\newblock In \emph{KDD workshop}, volume~10, pages 359--370. Seattle, WA, USA:,
  1994.

\bibitem[Boedihardjo et~al.(2016)Boedihardjo, Geng, Lyons, and Yang]{BGLY16}
Horatio Boedihardjo, Xi~Geng, Terry Lyons, and Danyu Yang.
\newblock The signature of a rough path: uniqueness.
\newblock \emph{Adv. Math.}, 293:\penalty0 720--737, 2016.

\bibitem[Brent(1973)]{Brent73}
Richard~P. Brent.
\newblock \emph{Algorithms for minimization without derivatives}.
\newblock Prentice-Hall, Inc., Englewood Cliffs, N.J., 1973.
\newblock Prentice-Hall Series in Automatic Computation.

\bibitem[Cass et~al.(2016)Cass, Driver, Lim, and Litterer]{CDLL16}
Thomas Cass, Bruce~K. Driver, Nengli Lim, and Christian Litterer.
\newblock On the integration of weakly geometric rough paths.
\newblock \emph{J. Math. Soc. Japan}, 68\penalty0 (4):\penalty0 1505--1524,
  2016.

\bibitem[Chan and Vasconcelos(2005)]{chan2005probabilistic}
Antoni~B Chan and Nuno Vasconcelos.
\newblock Probabilistic kernels for the classification of auto-regressive
  visual processes.
\newblock In \emph{2005 IEEE Computer Society Conference on Computer Vision and
  Pattern Recognition (CVPR'05)}, volume~1, pages 846--851. IEEE, 2005.

\bibitem[{Chevyrev} et~al.(2020){Chevyrev}, {Nanda}, and {Oberhauser}]{CNO20}
I.~{Chevyrev}, V.~{Nanda}, and H.~{Oberhauser}.
\newblock Persistence {P}aths and {S}ignature {F}eatures in {T}opological
  {D}ata {A}nalysis.
\newblock \emph{IEEE Transactions on Pattern Analysis and Machine
  Intelligence}, 42\penalty0 (1):\penalty0 192--202, 2020.

\bibitem[Chevyrev(2018)]{Chevyrev18}
Ilya Chevyrev.
\newblock Random walks and {L}\'evy processes as rough paths.
\newblock \emph{Probab. Theory Related Fields}, 170\penalty0 (3-4):\penalty0
  891--932, 2018.

\bibitem[Chevyrev and Friz(2019)]{CF19}
Ilya Chevyrev and Peter~K. Friz.
\newblock Canonical {RDE}s and general semimartingales as rough paths.
\newblock \emph{Ann. Probab.}, 47\penalty0 (1):\penalty0 420--463, 2019.

\bibitem[{Chevyrev} and {Kormilitzin}(2016)]{CK16}
Ilya {Chevyrev} and Andrey {Kormilitzin}.
\newblock {A Primer on the Signature Method in Machine Learning}.
\newblock \emph{ArXiv e-prints}, March 2016.

\bibitem[Chevyrev and Lyons(2016)]{ChevyrevLyons16}
Ilya Chevyrev and Terry Lyons.
\newblock {C}haracteristic {F}unctions of {M}easures on {G}eometric {R}ough
  {P}aths.
\newblock \emph{Ann. Probab.}, 44\penalty0 (6):\penalty0 4049--4082, 2016.

\bibitem[{Chevyrev} and {Ogrodnik}(2018)]{CO18}
Ilya {Chevyrev} and Marcel {Ogrodnik}.
\newblock A support and density theorem for {M}arkovian rough paths.
\newblock \emph{Electron. J. Probab.}, 23\penalty0 (56):\penalty0 16~pp., 2018.

\bibitem[Coutin and Qian(2002)]{coutin-qian-02}
Laure Coutin and Zhongmin Qian.
\newblock Stochastic analysis, rough path analysis and fractional {B}rownian
  motions.
\newblock \emph{Probab. Theory Related Fields}, 122\penalty0 (1):\penalty0
  108--140, 2002.
\newblock ISSN 0178-8051.

\bibitem[Cucker and Smale(2002)]{cucker2002mathematical}
Felipe Cucker and Steve Smale.
\newblock On the mathematical foundations of learning.
\newblock \emph{Bulletin of the American mathematical society}, 39\penalty0
  (1):\penalty0 1--49, 2002.

\bibitem[Cuturi and Doucet(2011)]{cuturi2011autoregressive}
Marco Cuturi and Arnaud Doucet.
\newblock Autoregressive kernels for time series.
\newblock \emph{arXiv preprint arXiv:1101.0673}, 2011.

\bibitem[Fawcett(2003)]{Fawcettthesis}
Thomas Fawcett.
\newblock \emph{Problems in stochastic analysis: connections between rough
  paths and non-commutative harmonic analysis}.
\newblock PhD thesis, University of Oxford, 2003.

\bibitem[Friedman and Rafsky(1979)]{friedman1979multivariate}
Jerome~H Friedman and Lawrence~C Rafsky.
\newblock Multivariate generalizations of the {W}ald-{W}olfowitz and {S}mirnov
  two-sample tests.
\newblock \emph{The Annals of Statistics}, pages 697--717, 1979.

\bibitem[Friz and Shekhar(2017)]{FrizShekhar17}
Peter Friz and Atul Shekhar.
\newblock {G}eneral rough integration, {L}{\'e}vy rough paths and a
  {L}{\'e}vy--{K}intchine type formula.
\newblock \emph{Ann. Probab.}, 45\penalty0 (4):\penalty0 2707–2765, 2017.

\bibitem[Friz and Victoir(2008)]{friz-victoir-subelliptic}
Peter Friz and Nicolas Victoir.
\newblock On uniformly subelliptic operators and stochastic area.
\newblock \emph{Probab. Theory Related Fields}, 142\penalty0 (3-4):\penalty0
  475--523, 2008.

\bibitem[Friz and Hairer(2014)]{friz2014course}
Peter~K. Friz and Martin Hairer.
\newblock \emph{A course on rough paths}.
\newblock Universitext. Springer, Cham, 2014.
\newblock With an introduction to regularity structures.

\bibitem[Friz and Victoir(2010)]{friz-victoir-book}
Peter~K. Friz and Nicolas~B. Victoir.
\newblock \emph{Multidimensional stochastic processes as rough paths}, volume
  120 of \emph{Cambridge Studies in Advanced Mathematics}.
\newblock Cambridge University Press, Cambridge, 2010.

\bibitem[Friz et~al.(2016)Friz, Gess, Gulisashvili, and Riedel]{FGGR16}
Peter~K. Friz, Benjamin Gess, Archil Gulisashvili, and Sebastian Riedel.
\newblock The {J}ain-{M}onrad criterion for rough paths and applications to
  random {F}ourier series and non-{M}arkovian {H}\"ormander theory.
\newblock \emph{Ann. Probab.}, 44\penalty0 (1):\penalty0 684--738, 2016.

\bibitem[Fukumizu et~al.(2009)Fukumizu, Gretton, Lanckriet, Sch{\"o}lkopf, and
  Sriperumbudur]{fukumizu2009kernel}
Kenji Fukumizu, Arthur Gretton, Gert~R Lanckriet, Bernhard Sch{\"o}lkopf, and
  Bharath~K Sriperumbudur.
\newblock Kernel choice and classifiability for rkhs embeddings of probability
  distributions.
\newblock In \emph{Advances in neural information processing systems}, pages
  1750--1758, 2009.

\bibitem[Giles(1971)]{giles1971generalization}
Robin Giles.
\newblock A generalization of the strict topology.
\newblock \emph{Transactions of the American Mathematical Society},
  161:\penalty0 467--474, 1971.

\bibitem[Gretton et~al.(2009)Gretton, Fukumizu, Harchaoui, and
  Sriperumbudur]{gretton2009fast}
Arthur Gretton, Kenji Fukumizu, Zaid Harchaoui, and Bharath~K Sriperumbudur.
\newblock A fast, consistent kernel two-sample test.
\newblock In \emph{Advances in neural information processing systems}, pages
  673--681, 2009.

\bibitem[Gretton et~al.(2012)Gretton, Borgwardt, Rasch, Sch{\"o}lkopf, and
  Smola]{gretton2012kernel}
Arthur Gretton, Karsten~M Borgwardt, Malte~J Rasch, Bernhard Sch{\"o}lkopf, and
  Alexander Smola.
\newblock A kernel two-sample test.
\newblock \emph{Journal of Machine Learning Research}, 13\penalty0
  (Mar):\penalty0 723--773, 2012.

\bibitem[Hampel et~al.(1986)Hampel, Ronchetti, Rousseeuw, and
  Stahel]{Hampel_robust}
Frank~R. Hampel, Elvezio~M. Ronchetti, Peter~J. Rousseeuw, and Werner~A.
  Stahel.
\newblock \emph{Robust statistics: the approach based on influence functions}.
\newblock Wiley Series in Probability and Mathematical Statistics: Probability
  and Mathematical Statistics. John Wiley \& Sons, Inc., New York, 1986.
\newblock ISBN 0-471-82921-8.

\bibitem[Huber and Ronchetti(2009)]{huber}
Peter~J. Huber and Elvezio~M. Ronchetti.
\newblock \emph{Robust statistics}.
\newblock Wiley Series in Probability and Statistics. John Wiley \& Sons, Inc.,
  Hoboken, NJ, second edition, 2009.
\newblock ISBN 978-0-470-12990-6.
\newblock \doi{10.1002/9780470434697}.
\newblock URL \url{https://doi.org/10.1002/9780470434697}.

\bibitem[Jones et~al.(2001--)Jones, Oliphant, and Peterson]{python}
Eric Jones, Travis Oliphant, and Pearu Peterson.
\newblock {SciPy}: Open source scientific tools for {Python}, 2001--.
\newblock URL \url{http://www.scipy.org/}.

\bibitem[Kiraly and Oberhauser(2019)]{OK16}
Franz~J. Kiraly and Harald Oberhauser.
\newblock Kernels for sequentially ordered data.
\newblock \emph{Journal of Machine Learning Research}, 20\penalty0
  (31):\penalty0 1--45, 2019.

\bibitem[Lang(2002)]{lang}
Serge Lang.
\newblock \emph{Algebra}.
\newblock Springer,, 2002.

\bibitem[Lehmann and Romano(2006)]{lehmann2006testing}
Erich~L Lehmann and Joseph~P Romano.
\newblock \emph{Testing statistical hypotheses}.
\newblock Springer Science \& Business Media, 2006.

\bibitem[Lyons(1998)]{lyons-98}
Terry Lyons.
\newblock Differential equations driven by rough signals.
\newblock \emph{Rev. Mat. Iberoamericana}, 14\penalty0 (2):\penalty0 215--310,
  1998.
\newblock ISSN 0213-2230.

\bibitem[Lyons and Ni(2015)]{ni2015expected}
Terry Lyons and Hao Ni.
\newblock Expected signature of {B}rownian motion up to the first exit time
  from a bounded domain.
\newblock \emph{Ann. Probab.}, 43\penalty0 (5):\penalty0 2729--2762, 2015.

\bibitem[Lyons and Qian(2002)]{lyons-qian-02}
Terry Lyons and Zhongmin Qian.
\newblock \emph{System control and rough paths}.
\newblock Oxford Mathematical Monographs. Oxford University Press, Oxford,
  2002.
\newblock Oxford Science Publications.

\bibitem[Lyons et~al.(2007)Lyons, Caruana, and L{\'e}vy]{lyons-04}
Terry Lyons, Michael Caruana, and Thierry L{\'e}vy.
\newblock \emph{Differential equations driven by rough paths}, volume 1908 of
  \emph{Lecture Notes in Mathematics}.
\newblock Springer, Berlin, 2007.

\bibitem[Moreno et~al.(2003)Moreno, Ho, and Vasconcelos]{moreno2003kullback}
Pedro~J Moreno, Purdy Ho, and Nuno Vasconcelos.
\newblock {A Kullback-Leibler Divergence Based Kernel for SVM Classification in
  Multimedia Applications.}
\newblock In \emph{NIPS}, pages 1385--1392, 2003.

\bibitem[Muandet et~al.(2017)Muandet, Fukumizu, Sriperumbudur, and
  Sch{\"o}lkopf]{muandet2017kernel}
Krikamol Muandet, Kenji Fukumizu, Bharath Sriperumbudur, and Bernhard
  Sch{\"o}lkopf.
\newblock Kernel mean embedding of distributions: A review and beyond.
\newblock \emph{Foundations and Trends{\textregistered} in Machine Learning},
  10\penalty0 (1-2):\penalty0 1--141, 2017.

\bibitem[M{\"u}ller(1997)]{muller1997integral}
Alfred M{\"u}ller.
\newblock Integral probability metrics and their generating classes of
  functions.
\newblock \emph{Advances in Applied Probability}, 29\penalty0 (2):\penalty0
  429--443, 1997.

\bibitem[Papavasiliou and Ladroue(2011)]{papavasiliou2011parameter}
Anastasia Papavasiliou and Christophe Ladroue.
\newblock Parameter estimation for rough differential equations.
\newblock \emph{The Annals of Statistics}, 39\penalty0 (4):\penalty0
  2047--2073, 2011.

\bibitem[Rachev(1991)]{rachev1991probability}
Svetlozar~Todorov Rachev.
\newblock \emph{Probability metrics and the stability of stochastic models},
  volume 269.
\newblock John Wiley \& Son Ltd, 1991.

\bibitem[Ramdas et~al.(2015)Ramdas, Reddi, P{\'o}czos, Singh, and
  Wasserman]{ramdas2015decreasing}
Aaditya Ramdas, Sashank~Jakkam Reddi, Barnab{\'a}s P{\'o}czos, Aarti Singh, and
  Larry~A Wasserman.
\newblock On the decreasing power of kernel and distance based nonparametric
  hypothesis tests in high dimensions.
\newblock In \emph{AAAI}, pages 3571--3577, 2015.

\bibitem[Reutenauer(1993)]{reutenauer-93}
Christophe Reutenauer.
\newblock \emph{Free {L}ie algebras}.
\newblock The Clarendon Press Oxford University Press, New York, 1993.
\newblock ISBN 0-19-853679-8.
\newblock Oxford Science Publications.

\bibitem[Sch{\"o}lkopf and Smola(2002)]{scholkopf2002learning}
Bernhard Sch{\"o}lkopf and Alexander~J Smola.
\newblock \emph{Learning with kernels: support vector machines, regularization,
  optimization, and beyond}.
\newblock MIT press, 2002.

\bibitem[Simon-Gabriel and Sch{\"o}lkopf(2018)]{SGS18}
Carl-Johann Simon-Gabriel and Bernhard Sch{\"o}lkopf.
\newblock Kernel distribution embeddings: Universal kernels, characteristic
  kernels and kernel metrics on distributions.
\newblock \emph{Journal of Machine Learning Research}, 19\penalty0
  (44):\penalty0 1--29, 2018.

\bibitem[{Simon-Gabriel} et~al.(2020){Simon-Gabriel}, {Barp}, and
  {Mackey}]{simongabriel2020metrizing}
Carl-Johann {Simon-Gabriel}, Alessandro {Barp}, and Lester {Mackey}.
\newblock {Metrizing Weak Convergence with Maximum Mean Discrepancies}.
\newblock \emph{arXiv e-prints}, art. arXiv:2006.09268, June 2020.

\bibitem[Sriperumbudur(2016)]{sriperumbudur2016optimal}
Bharath Sriperumbudur.
\newblock On the optimal estimation of probability measures in weak and strong
  topologies.
\newblock \emph{Bernoulli}, 22\penalty0 (3):\penalty0 1839--1893, 2016.

\bibitem[Sriperumbudur et~al.(2010)Sriperumbudur, Gretton, Fukumizu,
  Sch{\"o}lkopf, and Lanckriet]{sriperumbudur2010hilbert}
Bharath~K Sriperumbudur, Arthur Gretton, Kenji Fukumizu, Bernhard
  Sch{\"o}lkopf, and Gert~RG Lanckriet.
\newblock Hilbert space embeddings and metrics on probability measures.
\newblock \emph{Journal of Machine Learning Research}, 11\penalty0
  (Apr):\penalty0 1517--1561, 2010.

\bibitem[Sriperumbudur et~al.(2011)Sriperumbudur, Fukumizu, and
  Lanckriet]{sriperumbudur2011universality}
Bharath~K Sriperumbudur, Kenji Fukumizu, and Gert~RG Lanckriet.
\newblock Universality, characteristic kernels and rkhs embedding of measures.
\newblock \emph{Journal of Machine Learning Research}, 12\penalty0
  (Jul):\penalty0 2389--2410, 2011.

\bibitem[{Sutherland} et~al.(2016){Sutherland}, {Tung}, {Strathmann}, {De},
  {Ramdas}, {Smola}, and {Gretton}]{sutherland2016generative}
Dougal~J. {Sutherland}, Hsiao-Yu {Tung}, Heiko {Strathmann}, Soumyajit {De},
  Aaditya {Ramdas}, Alex {Smola}, and Arthur {Gretton}.
\newblock {Generative Models and Model Criticism via Optimized Maximum Mean
  Discrepancy}.
\newblock \emph{arXiv e-prints}, art. arXiv:1611.04488, November 2016.

\bibitem[Toth and Oberhauser(2020)]{2019arXiv190608215T}
Csaba Toth and Harald Oberhauser.
\newblock Bayesian {L}earning from {S}equential {D}ata using {G}aussian
  {P}rocesses with {S}ignature {C}ovariances.
\newblock In \emph{International Conference on Machine Learning}, pages
  9548--9560. PMLR, 2020.

\end{thebibliography}

\end{document}